\documentclass[11pt,a4paper]{article}
\usepackage[margin=28mm]{geometry}
\usepackage{amsmath,amssymb,amsthm,mathtools,microtype}
\usepackage[hidelinks]{hyperref}
\allowdisplaybreaks
\numberwithin{equation}{section}
\newtheorem{theorem}{Theorem}[section]
\newtheorem{lemma}{Lemma}[section]
\newtheorem{proposition}[lemma]{Proposition}
\theoremstyle{remark}
\newtheorem{remark}[lemma]{Remark}
\newcommand{\dd}{\,\mathrm dx}
\newcommand{\ds}{\,\mathrm ds}
\newcommand{\diver}{\operatorname{div}}
\newcommand{\curl}{\operatorname{curl}}
\newcommand{\Dt}{D_t}
\newcommand{\eps}{\varepsilon}
\newcommand{\Om}{\Omega}
\title{The draft of {\bf Global dynamics of large smooth solutions to the
two-dimensional compressible Navier-Stokes system with Navier-slip boundary conditions}\thanks{This work is partially supported by the National Natural Science Foundation of
China (12471227, 12494541).}}

\author{Xinyu, F{\small AN}$^{a}$,  Song J{\small IANG}$^{b}$, Jing L{\small I}$^{c,d}$, Boqiang L{\small \"U}$^{d} $  \thanks{email:  fanxinyu17@mails.ucas.ac.cn(X. Fan),  jiang@iapcm.ac.cn(S. Jiang), ajingli@gmail.com (J. Li), lvbq86@ncu.edu.cn (B.  L\"u) }  \\
{\normalsize a. School of Mathematics  Sciences, } \\
{\normalsize  Nankai University, Tianjin 300110, P. R. China;}
 \\ {\normalsize b. Institute of Applied Physics and Computational Mathematics, }\\\ {\normalsize Beijing 100088, P. R. China;}\\
{\normalsize c. Institute of Applied Mathematics \& State Key Laboratory of Mathematical Sciences,  AMSS,} \\ {\normalsize \&   Hua Loo-Keng Key Laboratory of Mathematics,}\\
{\normalsize  Chinese Academy of Sciences,    Beijing 100190, P. R. China};\\
{\normalsize d. School of Mathematics and Computer Sciences,}\\  {\normalsize  Nanchang University, Nanchang 330031, P. R. China.}
}
\date{}
\hypersetup{
  pdftitle={Global strong solutions and large-time behavior for the two-dimensional compressible Navier-Stokes equations with Navier slip},
  pdfauthor={Jing Li;  },
  pdfkeywords={compressible Navier-Stokes equations; Navier slip boundary conditions; density-dependent bulk viscosity; global strong solutions; exponential decay; vacuum}
}
\begin{document}
\maketitle
\begin{abstract}
The global existence of a unique strong solution to the two-dimensional
barotropic compressible Navier--Stokes equations with Navier slip
boundary conditions is established in a smooth bounded domain with large initial data and vacuum.
The bulk viscosity $\lambda$ satisfies a general structural condition, which includes the case $\lambda=\rho^\beta$ with $\beta>1$. 
Exponential decay of the density and velocity is proved in every
finite Lebesgue and first-order Sobolev space, respectively.
For strictly positive initial density, uniform positivity and
exponential decay of density gradients and second velocity derivatives
are established. In contrast, if
an initial vacuum occurs, then the vacuum persists and causes exponential growth of
density gradients.
\end{abstract}

\section{Introduction and main result}

Let $\Om\subset\mathbb R^2$ be a bounded connected smooth domain,
possibly with holes. We consider the compressible Navier--Stokes equations
\begin{equation}\label{eq:system}
\begin{cases}
 \rho_t+\diver(\rho u)=0,\\
 (\rho u)_t+\diver(\rho u\otimes u)+\nabla P
     =\mu\Delta u+\nabla((\mu+\lambda(\rho))\diver u)
\end{cases}
\end{equation}
in $\Om\times(0,\infty)$, where $P=\rho^\gamma$, $\mu>0$ and
$\gamma>1$. Write $\dot u=u_t+u\cdot\nabla u$.
Let $n$ be the outward unit normal, and set
\[
 J=\begin{pmatrix}0&1\\-1&0\end{pmatrix},\qquad
 \tau=Jn,\qquad u_\tau=u\cdot\tau,\qquad
 \omega=\partial_2u_1-\partial_1u_2.
\]
The classical Navier slip law~\cite{Navier1827} is
\begin{equation}\label{eq:navier-stress}
 u\cdot n=0,\qquad 2\mu D(u)n\cdot\tau+f u_\tau=0
 \quad\hbox{on }\partial\Om\times(0,\infty),
\end{equation}
where $D(u)=(\nabla u+\nabla u^{\mathsf T})/2$ and $f$ is the
smooth boundary friction coefficient. Define the signed curvature by
$\partial_\tau n=\kappa\tau$. On circular boundary components of
radius $R$, this gives $\kappa=1/R$ for an outer boundary and
$\kappa=-1/R$ for an inner boundary. Since $u\cdot n=0$,
\[
\begin{aligned}
 (\partial_\tau u)\cdot n&=-\kappa u_\tau,\\
 2D(u)n\cdot\tau
 &=(\partial_nu)\cdot\tau+(\partial_\tau u)\cdot n
   =\omega-2\kappa u_\tau.
\end{aligned}
\]
Thus \eqref{eq:navier-stress} is equivalent to the vorticity form
(see Coron~\cite[p.~36]{Coron1996})
\begin{equation}\label{eq:slip}
 u\cdot n=0,\qquad \omega=-\alpha u_\tau
 \quad\hbox{on }\partial\Om\times(0,\infty),
\end{equation}
where
\[
 \alpha=\frac f\mu-2\kappa.
\]
We assume $\alpha\ge0$, as in Fan, Li, and Li~\cite{FanLiLi2022}.
The hypotheses $\alpha\ge0$ and $\alpha\not\equiv0$ in the main
theorem therefore mean
\begin{equation}\label{eq:friction-condition}
 f\ge2\mu\kappa\quad\hbox{on }\partial\Om,\qquad
 f-2\mu\kappa\not\equiv0.
\end{equation}
Thus $\alpha$ is the coefficient after the curvature correction.
If nonnegative physical friction $f\ge0$ is also required, the pointwise
condition is $f\ge\max\{0,2\mu\kappa\}$, with the same nonzero-difference
condition. In particular, zero physical friction $f=0$ corresponds
to $\alpha=-2\kappa$.

The initial conditions are
\begin{equation}\label{eq:initial}
 \rho|_{t=0}=\rho_0,\qquad (\rho u)|_{t=0}=\rho_0u_0.
\end{equation}
All norms are on $\Om$, and $\bar g=|\Om|^{-1}\int_\Om g\dd$.

Assume that
\begin{equation}\label{eq:coefficient}
 \lambda\in C^1([0,\infty))\cap C^2((0,\infty)),\quad
 \lambda\ge0,\quad
 |r\lambda'(r)|\le C_\lambda(1+\lambda(r)),\quad r>0.
\end{equation}
Define
\[
 \theta(r)=2\mu\log r+\int_1^r\frac{\lambda(z)}z\,\mathrm dz,
 \qquad r>0.
\]
For some $a\in(1/2,1]$, $c_\lambda>0$ and $r_\lambda\ge2$, assume
\begin{equation}\label{eq:growth}
 \theta(r)\ge c_\lambda r[\log(e+r)]^a\quad(r\ge r_\lambda).
\end{equation}
Our main result is stated as follows.

\begin{theorem}\label{thm:strong}
Assume \eqref{eq:coefficient}--\eqref{eq:growth} and
\eqref{eq:friction-condition}. Let $q>2$,
$0\le\rho_0\in W^{1,q}(\Om)$ with $M=\int_\Om\rho_0\dd>0$, and
$u_0\in H^1(\Om)$ with $u_0\cdot n=0$. Then the problem
\eqref{eq:system}--\eqref{eq:initial} admits a unique global strong
solution. For every finite $T$,
\begin{equation}\label{str:class}
\begin{aligned}
 &\rho\in C([0,T];W^{1,q}),\qquad \rho_t\in L^\infty(0,T;L^2),\\
 &u\in L^\infty(0,T;H^1)\cap L^2(0,T;H^2)
                \cap L^{(q+1)/q}(0,T;W^{2,q}),\\
 &\sqrt t\,u\in L^2(0,T;W^{2,q}),\qquad
    \sqrt t\,u_t\in L^2(0,T;H^1),\\
 &\sqrt\rho\,u_t\in L^2(0,T;L^2),\qquad
     \rho u,\sqrt\rho\,u\in C([0,T];L^2).
\end{aligned}
\end{equation}
The initial density and momentum are attained in these continuous
spaces, and \eqref{eq:slip} holds in the trace sense for almost every
$t>0$.
Moreover, with $\sigma(t)=\min\{1,t\}$,
\begin{equation}\label{str:uniform}
\begin{aligned}
 &\operatorname*{ess\,sup}_{t>0}
       \bigl(\|\rho(t)\|_\infty+\|u(t)\|_{H^1}^2\bigr)
       +\int_0^\infty\|\sqrt\rho\,\dot u\|_2^2\,\mathrm dt\\
 &\quad+\operatorname*{ess\,sup}_{t>0}
       \sigma(t)\|\sqrt\rho\,\dot u(t)\|_2^2
       +\int_0^\infty\sigma(t)\|\nabla\dot u\|_2^2\,\mathrm dt\le C.
\end{aligned}
\end{equation}
The constant depends only on $\Om,\alpha,\mu,\gamma,M$,
$\|\rho_0\|_\infty$, $\|u_0\|_{H^1}$, the constants in
\eqref{eq:coefficient}--\eqref{eq:growth}, and a bound for $\lambda$
on $[0,2]$.
\end{theorem}

\begin{theorem}[Large-time behavior]\label{thm:longtime}
Under the hypotheses of Theorem~\ref{thm:strong}, set
\[
 \bar\rho=\frac{M}{|\Om|},\qquad
 G=(2\mu+\lambda(\rho))\diver u-(P-\bar P).
\]
For every $1\le p<\infty$, there are $C_p,c_p>0$ such that
\begin{equation}\label{lt:relaxation}
 \|\rho(t)-\bar\rho\|_p+\|u(t)\|_{W^{1,p}}
       \le C_pe^{-c_pt},\qquad t\ge1.
\end{equation}
Moreover, for almost every $t\ge1$,
\begin{equation}\label{lt:acceleration}
\begin{aligned}
 &\|\sqrt\rho\,\dot u(t)\|_2+\|\sqrt\rho\,u_t(t)\|_2
       +\|G(t)\|_{H^1}+\|\omega(t)\|_{H^1}\\
 &\qquad+\left(\int_t^\infty\|\dot u(s)\|_{H^1}^2\,\mathrm ds\right)^{1/2}
       \le Ce^{-ct}.
\end{aligned}
\end{equation}
These constants have the low-order dependence of
Theorem~\ref{thm:strong}, with an additional dependence on $p$ in
\eqref{lt:relaxation}.

If $\min_{\overline\Om}\rho_0>0$, then
\begin{equation}\label{lt:positive}
\begin{aligned}
 &\inf_{t\ge0,\,x\in\overline\Om}\rho(t,x)>0,\\
 &\|\rho(t)-\bar\rho\|_{W^{1,q}}+\|\rho_t(t)\|_q
        +\|u(t)\|_{H^2}+\|u_t(t)\|_2\le Ce^{-ct},\\
 &\int_t^\infty\left(\|u(s)\|_{W^{2,q}}^2
             +\|u_t(s)\|_{H^1}^2+\|\nabla u(s)\|_\infty\right)\,\mathrm ds
             \le Ce^{-ct}.
\end{aligned}
\end{equation}
The pointwise derivative estimates hold for almost every $t\ge1$;
the density estimate and the tail estimate hold for every $t\ge1$.
Here the constants may additionally depend on $q$,
$\|\rho_0\|_{W^{1,q}}$, and $(\min\rho_0)^{-1}$.

If instead $\rho_0(x_0)=0$ for some $x_0\in\overline\Om$, then
\begin{equation}\label{lt:vacuum}
 \min_{\overline\Om}\rho(t)=0,\qquad
 \|\rho(t)-\bar\rho\|_\infty\ge\bar\rho,\qquad t\ge0.
\end{equation}
For every $2<r\le q$, there are $C_r,c_r>0$ such that
\begin{equation}\label{lt:gradient-growth}
 \|\nabla\rho(t)\|_r\ge C_r^{-1}e^{c_rt},\qquad t\ge0.
\end{equation}
For every $0\le\varepsilon<1$, the low-density set satisfies
\begin{equation}\label{lt:vacuum-set}
 |\{x\in\Om:\rho(t,x)\le\varepsilon\bar\rho\}|
       \le C_\varepsilon e^{-ct},\qquad t\ge0.
\end{equation}
\end{theorem}

\begin{remark}[Comparison with earlier results]\label{rem:comparison}
For $\lambda(\rho)=\rho^\beta$, Fan, Li, and Li~\cite[Theorem~1.1]{FanLiLi2022}
proved global strong solvability in simply connected bounded domains
for $\beta>4/3$, with large initial data and vacuum allowed.
Fan, Jiang, and Li~\cite[Theorem~1.2]{FanJiangLi2024} treated
multiply connected domains under the same power restriction and
obtained uniform density bounds and exponential decay when
$\alpha\ge0$ and $\alpha\not\equiv0$.
Theorem~\ref{thm:strong} extends the strong existence range in these
works to every $\beta>1$ and allows general bulk viscosities satisfying
\eqref{eq:coefficient}--\eqref{eq:growth}, together with the uniform
estimates \eqref{str:uniform} and the large-time conclusions of
Theorem~\ref{thm:longtime}.
\end{remark}

\begin{remark}[Admissible bulk viscosities]\label{rem:viscosity}
Conditions \eqref{eq:coefficient}--\eqref{eq:growth} include
$\lambda(r)=r^\beta$ for every $\beta>1$ and
$\lambda(r)=r[\log(e+r)]^\eta$ for every $\eta>1/2$.
The latter grows more slowly than every superlinear power.
Neither monotonicity of $\lambda$ nor a uniform bound for $\lambda''$
near vacuum is required.
\end{remark}

\begin{remark}[Method and domain topology]\label{rem:method}
The proofs in \cite{FanLiLi2022,FanJiangLi2024} use Green function
representations and conformal mappings. Here the effective viscous
flux is estimated through the variational Neumann operator and the
transport commutator in Lemma~\ref{mom:commutator}.
The local elliptic estimates in Lemma~\ref{mom:potential} require no
estimates for the domain Green function.
The Hodge decomposition retains the finite-dimensional harmonic part,
so holes are handled without any symmetry assumption.
The uniform bounds depend only on the low-order data listed in
Theorem~\ref{thm:strong}, independently of a positive density lower bound.
\end{remark}

\begin{remark}[Boundary dissipation]\label{rem:general-slip}
Under $\alpha\ge0$, positivity on one open boundary arc suffices for
the coercivity \eqref{uni:coercive}; strict positivity on the whole
boundary is unnecessary.
For an arbitrary smooth $\alpha\ge0$, equivalently $f\ge2\mu\kappa$,
including $\alpha\equiv0$ (that is, $f=2\mu\kappa$),
the problem still has a unique global strong solution satisfying
\eqref{str:class}. The low-order estimates on $[0,T]$ are then those
of Proposition~\ref{thm:main}, with a constant $C_T$.
The time-uniform estimate \eqref{str:uniform} is not asserted when
$\alpha\equiv0$, since the dissipation may have a nontrivial harmonic
kernel on a multiply connected domain.
Nontrivial stationary circulations occur on annuli when
$\alpha\equiv0$; see \cite[Theorem~1.1(b)]{FanJiangLi2024}.
Thus decay to $(\bar\rho,0)$ is not asserted in this case.
\end{remark}

\begin{remark}[Cases not covered]\label{rem:remaining}
The endpoint $\lambda(r)=r[\log(e+r)]^{1/2}$ and the linear law
$\lambda(r)=r$ are not covered by \eqref{eq:growth}.
The strict condition $a>1/2$ enters the moment estimate through the
choice $(2a-1)k\ge1$ for a finite $k$; its sharpness is not established
here.
The present argument also leaves untreated nonnegative physical
friction $f\ge0$ without the restriction $f\ge2\mu\kappa$.
The corresponding curvature-corrected coefficient
$\alpha=f/\mu-2\kappa$ may be negative, which lies outside the present
assumptions.
\end{remark}

\section{A priori estimates}

For a smooth positive solution, set
\[
 F=(2\mu+\lambda(\rho))\diver u-P,\qquad
 \nabla^\perp=(\partial_2,-\partial_1).
\]
The momentum equation becomes
\begin{equation}\label{eq:flux}
 \rho\dot u=\nabla F+\mu\nabla^\perp\omega.
\end{equation}

Define
\begin{equation}\label{eq:energy}
\begin{aligned}
 A_3^2&=\int_\Om\bigl((2\mu+\lambda)(\diver u)^2+\mu\omega^2\bigr)\dd
            +\mu\int_{\partial\Om}\alpha u_\tau^2\ds,\\
 y=A_2^2&=\int_\Om\rho|\dot u|^2\dd,\\
 E&=e+\int_\Om\left(\mu\omega^2+\frac{F^2}{2\mu+\lambda}\right)\dd
            +\mu\int_{\partial\Om}\alpha u_\tau^2\ds,\\
 R_t&=\max\{2,\sup_{0\le s\le t}\|\rho(s)\|_\infty\}.
\end{aligned}
\end{equation}

\begin{proposition}\label{thm:main}
Fix $T>0$, $C_0\ge1$, and a positive total mass $M$.
Let $(\rho,u)$ be a smooth solution of
\eqref{eq:system}--\eqref{eq:slip} on $[0,S]$, $0<S\le T$, with
$\rho>0$ and
\[
 \int_\Om\rho_0\dd=M,\qquad
 \|\rho_0\|_\infty+\|u_0\|_{H^1}\le C_0.
\]
Under \eqref{eq:coefficient}--\eqref{eq:growth},
\begin{equation}\label{eq:main}
\begin{aligned}
 &\sup_{0\le t\le S}
   \bigl(\|\rho(t)\|_\infty+\|u(t)\|_{H^1}^2+E(t)+A_3^2(t)\bigr)
   +\int_0^S A_2^2\,\mathrm dt\\
 &\qquad+\sup_{0<t\le S}tA_2^2(t)
       +\int_0^S t\|\nabla\dot u\|_2^2\,\mathrm dt\le C_T.
\end{aligned}
\end{equation}
The constant depends on the fixed domain, boundary coefficient, mass,
$\mu,\gamma,C_0,T$, the constants in
\eqref{eq:coefficient}--\eqref{eq:growth}, and a fixed bound for
$\lambda$ on $[0,2]$. It is independent of $S$, the positive density
lower bound, and higher initial norms.
\end{proposition}

Since $\theta'(r)=(2\mu+\lambda(r))/r>0$, condition
\eqref{eq:coefficient} implies
\begin{equation}\label{eq:coeffconsequence}
 1+\lambda(r)\le C(1+r^b),\qquad
 \lambda(r)\le C(1+\theta(r)_+),\qquad b=\max\{1,C_\lambda\}.
\end{equation}
The first bound follows by integration. For the second, integrate
$\lambda(z)/z$ over $[r/2,r]$ for $r\ge2$ and use
$1+\lambda(z)\ge2^{-C_\lambda}(1+\lambda(r))$.

Multiplication of momentum by $u$ gives the basic energy identity
including its boundary dissipation:
\begin{equation}\label{eq:basic}
 \frac{\mathrm d}{\mathrm dt}\int_\Om
       \left(\frac12\rho|u|^2+\frac{\rho^\gamma}{\gamma-1}\right)\dd
       +A_3^2=0,\qquad
 \|u\|_{H^1}\le C(1+A_3).
\end{equation}
The positive mass, basic density energy, and weighted div--curl
inequality give the last bound, including the harmonic velocities.
Transport boundary terms vanish since $u\cdot n=0$; the viscous
boundary term equals $-\mu\int_{\partial\Om}\alpha u_\tau^2\ds$.
All later constants $C_T$ have the dependence in
Proposition~\ref{thm:main}; $C_{\tau,T}$ may additionally depend on $\tau>0$.

\subsection{Elliptic estimates}
\label{sec:elliptic}

For a vector field $f$, define the mean-zero scalar $Qf$ by
\begin{equation}\label{mom:Q}
 \int_\Om\nabla Qf\cdot\nabla\phi\dd
       =-\int_\Om f\cdot\nabla\phi\dd,\qquad \overline{Qf}=0.
\end{equation}
Thus $Q\nabla g=-g+\bar g$.

\begin{lemma}[Quantitative elliptic estimates]\label{mom:potential}
On the fixed smooth domain, for every $p\ge2$,
\begin{equation}\label{mom:elliptic}
 \|\nabla Qf\|_p\le Cp\|f\|_p.
\end{equation}
The same estimate holds for the gradient of the zero-Dirichlet
solution with divergence data. For the scalar homogeneous Dirichlet
or mean-zero Neumann Poisson problem $-\Delta w=g$,
\begin{equation}\label{ell:Hessian}
 \|D^2w\|_p\le Cp^6\|g\|_p .
\end{equation}
The Neumann equation has the usual mean-zero compatibility condition.
Consequently, tangent vector fields satisfy
\begin{equation}\label{mom:Hodge}
 \|\nabla v\|_p\le Cp^6(\|\diver v\|_p+\|\curl v\|_p)
                          +C\|v\|_2.
\end{equation}
For each fixed $1<q<2$, one also has
\begin{equation}\label{mom:potential-low}
 \|Qf\|_{2q/(2-q)}\le C_q\|f\|_q.
\end{equation}
\end{lemma}

\begin{proof}
The operator $\nabla Q$ is self-adjoint and has $L^2$ norm at most
one. To prove its weak $(1,1)$ bound, consider a
mean-zero vector bad block $b$ supported in $B_r\cap\Om$ and an
$L^2$ vector field $g$ supported in an annulus at distance comparable
to $d\ge8r$.

Near $B_r$, $Qg$ is harmonic with homogeneous Neumann boundary data.
The interior and boundary derivative estimates give
\[
 \operatorname{osc}_{B_r\cap\Om}\nabla Qg
 \le Crd^{-2}\|\nabla Qg\|_{L^2(B_{cd}\cap\Om)}
 \le Crd^{-2}\|g\|_2.
\]
For $d$ above a fixed chart radius, a finite cover by fixed-scale
homogeneous charts separated from $\operatorname{supp}g$ gives
oscillation at most $Cr\|g\|_2$ for small $r$, and $C\|g\|_2$ otherwise.
Since $d\le\operatorname{diam}\Om$ and $r$ is bounded below in the
latter case, the same displayed bound follows.
Self-adjointness and the cancellation of $b$ yield
\[
 \|\nabla Qb\|_{L^2(\text{annulus }d)}
       \le Crd^{-2}\|b\|_1,
 \qquad
 \|\nabla Qb\|_{L^1(\text{annulus }d)}\le C(r/d)\|b\|_1.
\]
Summing the dyadic annuli and applying the Calderon--Zygmund
decomposition on the doubling space $\overline\Om$, with the $L^2$
estimate for the good part, proves weak $(1,1)$.
Interpolation near exponent one, duality, and the $L^2$ bound give
\eqref{mom:elliptic}. The Dirichlet gradient projection is treated
identically, using homogeneous local Dirichlet estimates.

For \eqref{ell:Hessian}, choose normal coordinates near the boundary.
Even reflection for homogeneous Neumann data and odd reflection for
homogeneous Dirichlet data yield an elliptic equation with Lipschitz
principal coefficients and bounded lower order coefficients. In two
dimensions the normal coordinates have no tangential--normal cross
coefficient. A frozen constant-coefficient Fourier estimate has norm
at most $Cp$. Localize at radius $c/p$, with $c$ sufficiently small:
the principal coefficient oscillation is $O(c/p)$ and is absorbed.
Summing a cover of bounded
overlap gives
\[
 \|D^2w\|_p\le Cp\|g\|_p+Cp^2\|\nabla w\|_p+Cp^3\|w\|_p.
\]
Uniform interpolation
$\|\nabla w\|_p\le\delta\|D^2w\|_p+C\delta^{-1}\|w\|_p$
with $\delta\asymp p^{-2}$ absorbs the gradient term.
The fixed $H^2\hookrightarrow L^\infty$ estimate gives
$\|w\|_p\le C\|g\|_2\le C\|g\|_p$ uniformly for $p\ge2$,
proving \eqref{ell:Hessian}.

For \eqref{mom:Hodge}, solve
\[
 \Delta\phi=\diver v,\quad\partial_n\phi=0,\quad\bar\phi=0,
 \qquad \Delta\chi=\curl v,\quad\chi|_{\partial\Om}=0.
\]
The remainder $v-\nabla\phi-\nabla^\perp\chi$ is the orthogonal
projection onto the finite-dimensional space of smooth tangent
harmonic fields. Its gradient has $L^p$ norm at most $C\|v\|_2$,
uniformly for $p\ge2$. Apply \eqref{ell:Hessian} to both potentials
to obtain \eqref{mom:Hodge}. Fixed-exponent Neumann regularity and
Sobolev embedding give \eqref{mom:potential-low}.
\end{proof}

If $\rho\ge0$, $\int_\Om\rho\dd=M>0$, and
$\|\rho\|_\gamma\le C_0$ for $\gamma>1$, then
\begin{equation}\label{ell:weightedHodge}
 \|v\|_{H^1}\le C\bigl(\|\diver v\|_2+\|\curl v\|_2
                                    +\|\sqrt\rho\,v\|_2\bigr),
 \qquad v\cdot n=0.
\end{equation}
The constant depends only on $\Om,M,C_0,\gamma$. Indeed, a normalized
counterexample sequence converges strongly in $L^{2\gamma/(\gamma-1)}$,
while the densities converge weakly in $L^\gamma$ with mass $M$.
By \eqref{mom:Hodge}, the velocity limit is a nonzero tangent harmonic
field. The vanishing weighted norm forces it to vanish where the
limiting density is positive, a set of positive measure, contradicting
analyticity.

If $\alpha\ge0$ is smooth and not identically zero, then
\begin{equation}\label{uni:coercive}
 \|v\|_{H^1}^2\le C\left(\|\diver v\|_2^2+\|\curl v\|_2^2
                 +\int_{\partial\Om}\alpha v_\tau^2\ds\right),
 \qquad v\cdot n=0,
\end{equation}
with $C$ depending only on $\Om,\alpha$; see
\cite[Lemma~3.4]{FanJiangLi2024}.
Indeed, a normalized counterexample sequence and the Hodge decomposition
would give a nonzero tangent harmonic limit $h$ with $h_\tau=0$ on
an open boundary arc where $\alpha>0$. Locally write
$h=\nabla^\perp\chi$. Then $\chi$ is harmonic, constant on the arc,
and has zero normal derivative there. Extending $\chi$ minus that
constant by zero across the arc gives a harmonic function vanishing
on an open set. Unique continuation forces $h=0$, a contradiction.

\begin{lemma}[Flux and vorticity gradients]\label{mom:fluxelliptic-lemma}
For every $p\ge2$,
\begin{equation}\label{mom:fluxelliptic}
 \|\nabla F\|_p+\|\nabla\omega\|_p
 \le Cp\bigl(\|\rho\dot u\|_p+\|u\|_{W^{1,p}}\bigr).
\end{equation}
\end{lemma}
\begin{proof}
Extend $\alpha,\tau$ smoothly into $\Om$. Then $\omega+\alpha u\cdot\tau$
has zero trace and satisfies
\[
 \mu\Delta(\omega+\alpha u\cdot\tau)
     =\curl(\rho\dot u)+\mu\Delta(\alpha u\cdot\tau).
\]
Since $\|\alpha u\cdot\tau\|_{W^{1,p}}\le C\|u\|_{W^{1,p}}$
uniformly in $p$, the
Dirichlet divergence-data estimate in Lemma~\ref{mom:potential}
gives the vorticity bound; the momentum identity gives the flux bound.
\end{proof}

\begin{lemma}[Transport commutator]
\label{mom:commutator}
Define
\[
 \mathfrak C(u,f)=u\cdot\nabla Qf-Q\diver(f\otimes u),
 \qquad (\diver(f\otimes u))_j=\partial_i(f_j u_i).
\]
For a tangent $u$ and every smooth test function,
\begin{equation}\label{ell:commidentity}
\begin{aligned}
 &\int_\Om\nabla\mathfrak C(u,f)\cdot\nabla\phi\dd\\
 &=\int_\Om\left[(\nabla u+\nabla u^{\mathsf T}-(\diver u)I)\nabla Qf
                            +(\nabla u)f\right]\cdot\nabla\phi\dd.
\end{aligned}
\end{equation}
For each fixed $m>2$,
\begin{equation}\label{mom:commutator-bound}
 \|\mathfrak C(u,f)\|_m\le C_m\|\nabla u\|_2\|f\|_m.
\end{equation}
\end{lemma}
\begin{proof}
The field $f+\nabla Qf$ has zero divergence and normal trace by
\eqref{mom:Q}. Integrate the transport terms along $u$, using
$u\cdot n=0$, and then along $f+\nabla Qf$. This gives
\[
\begin{aligned}
 &\int_\Om\nabla\mathfrak C(u,f)\cdot\nabla\phi\dd\\
 &=\int_\Om(\partial_k u_i\partial_i Qf-(\diver u)\partial_k Qf)\partial_k\phi\dd
          -\int_\Om u_i(f_j+\partial_j Qf)\partial_{ij}\phi\dd\\
 &=\int_\Om(\partial_k u_i\partial_i Qf-(\diver u)\partial_k Qf)\partial_k\phi\dd
          +\int_\Om\partial_j u_i(f_j+\partial_j Qf)\partial_i\phi\dd.
\end{aligned}
\]
Neumann regularity and H\"older's inequality give
$\|\nabla\mathfrak C(u,f)\|_r\le C_m\|\nabla u\|_2\|f\|_m$ for
$1/r=1/2+1/m$. Its mean is
$\overline{\mathfrak C(u,f)}=-|\Om|^{-1}\int_\Om(\diver u)Qf\dd$
and has the same bound.
The embedding $W^{1,r}\hookrightarrow L^m$ completes the proof.
\end{proof}

\subsection{Density moments}\label{mom:section}

\begin{lemma}[The harmonic boundary contribution]\label{mom:boundary}
Let $B=Q(\nabla^\perp\omega)$. Then, for every finite $p\ge1$,
\begin{equation}\label{mom:boundary-bound}
 \|B\|_{H^1}+\|B\|_p\le C_{\Om,\alpha,p}\|u\|_{H^1}.
\end{equation}
\end{lemma}
\begin{proof}
Since $\diver\nabla^\perp\omega=0$ and
$n\cdot\nabla^\perp\omega=-\partial_\tau\omega$,
\eqref{mom:Q} and the slip condition give
\[
 \Delta B=0,\qquad
 \partial_nB=\partial_\tau\omega
            =-\partial_\tau(\alpha u_\tau),\qquad \bar B=0.
\]
Since $\int_{\partial\Om}\partial_\tau(\alpha u_\tau)\ds=0$,
Neumann regularity and the trace theorem give
\[
 \|B\|_{H^1}
  \le C\|\partial_\tau(\alpha u_\tau)\|_{H^{-1/2}(\partial\Om)}
  \le C\|\alpha u_\tau\|_{H^{1/2}(\partial\Om)}
  \le C\|u\|_{H^1}.
\]
Sobolev embedding gives the $L^p$ bound.
\end{proof}

Fix an integer $k\ge2$ such that
\begin{equation}\label{mom:k}
 k+1\ge4\gamma,\qquad (2a-1)k\ge1,
\end{equation}
where $a>1/2$ is the exponent in \eqref{eq:growth}.

\begin{lemma}[A fixed Orlicz density moment]\label{mom:moment}
For the fixed $k$ in \eqref{mom:k},
\begin{equation}\label{mom:moment-bound}
 \sup_{0\le t\le S}\int_\Om\rho\theta(\rho)_+^k\dd
 +\sup_{0\le t\le S}\int_\Om
       \rho^{k+1}[\log(e+\rho)]^{ak}\dd\le C_T.
\end{equation}
The constant has the dependence in Proposition~\ref{thm:main};
it is independent of the density lower bound and of $R$.
\end{lemma}
\begin{proof}
Put $\psi=Q(\rho u)$.
Applying $Q$ to conservative momentum gives
\[
 \Dt\psi=-F+\bar F+\mu B+\mathfrak C(u,\rho u).
\]
The continuity equation gives $\Dt\theta(\rho)+P=-F$. Hence
\begin{equation}\label{mom:transport}
 \Dt(\theta(\rho)-\psi)+P
    =-\bar F-\mu B-\mathfrak C(u,\rho u).
\end{equation}
Set
\[
 \mathcal M_k(t)=\int_\Om\rho(\theta(\rho)-\psi)_+^k\dd,
 \qquad m=k+1.
\]
The growth condition implies
\[
 \int_\Om\rho^m[\log(e+\rho)]^{ak}\dd
       \le C\left(1+\int_\Om\rho\theta(\rho)_+^k\dd\right).
\]
For $q=2m/(m+1)$, the potential estimate and the basic kinetic energy
bound yield
\[
 \|\psi\|_{2m}\le C_m\|\rho u\|_q
       \le C_m\|\rho\|_m^{1/2}.
\]
Interpolation with the fixed mass gives
$\|\rho\|_q\le C_M\|\rho\|_m^{1/2}$. Consequently,
\[
 \int_\Om\rho|\psi|^k\dd
  \le\|\rho\|_q\|\psi\|_{2m}^k
  \le C_m\|\rho\|_m^{m/2}
  \le C_m\left(1+\int_\Om\rho\theta(\rho)_+^k\dd\right)^{1/2}.
\]
Using $\theta_+\le(\theta-\psi)_++|\psi|$ and absorbing this square
root gives
\begin{equation}\label{mom:comparison}
 \int_\Om\rho\theta(\rho)_+^k\dd
 +\int_\Om\rho^m[\log(e+\rho)]^{ak}\dd
       \le C_m(1+\mathcal M_k).
\end{equation}

By \eqref{eq:basic} and the Trudinger inequality,
\[
 \int_\Om\exp\!\left(\frac{c|u|^2}{1+A_3^2}\right)\dd\le C.
\]
Apply the Orlicz--Young inequality
\[
 xy\le C_mx[\log(e+x)]^{m/2}
               +C_m\exp(c_my^{2/m})
\]
with $x=\rho^m$, $y=(c_*|u|/(1+A_3))^m$, and $c_*$ small.
Since $ak\ge m/2$ by \eqref{mom:k}, \eqref{mom:comparison} gives
\begin{equation}\label{mom:momentum}
 \|\rho u\|_m\le C_m(1+A_3)(1+\mathcal M_k)^{1/m}.
\end{equation}
Lemmas~\ref{mom:commutator} and \ref{mom:boundary} give
\[
 \|\mathfrak C(u,\rho u)\|_m
    \le C_m (1+A_3)^2(1+\mathcal M_k)^{1/m},\qquad
 \|B\|_m\le C_m(1+A_3).
\]

Test \eqref{mom:transport} by $\rho(\theta-\psi)_+^{k-1}$.
Since $u\cdot n=0$ and $P\ge0$, discard the pressure term and use
\[
 \int_\Om\rho(\theta-\psi)_+^{k-1}|Z|\dd
       \le\mathcal M_k^{1-1/k}\|\rho\|_m^{1/k}\|Z\|_m.
\]
Since $1-1/k+1/[k(k+1)]+1/(k+1)=1$, the commutator contribution
is bounded by $C_m(1+A_3^2)(1+\mathcal M_k)$.
The $B$ term has exponent $k/(k+1)<1$ and satisfies the same bound.

For the mean, since $\int_\Om\diver u\dd=0$,
\[
 \bar F=|\Om|^{-1}\int_\Om\lambda\diver u\dd-\bar P.
\]
By \eqref{eq:coeffconsequence}, $\theta_+=0$ for $\rho\le1$, and
\eqref{mom:comparison},
\[
 \int_\Om\lambda(\rho)\dd
   \le C\left(1+\left[\int_\Om\rho\theta(\rho)_+^k\dd\right]^{1/k}\right),
 \qquad
 |\bar F|\le C+C_mA_3(1+\mathcal M_k)^{1/(2k)}.
\]
By H\"older's inequality in the mass measure, the mean contribution
has exponent $1-1/(2k)<1$. Hence
\[
 \mathcal M_k'\le C_m(1+A_3^2)(1+\mathcal M_k).
\]
The initial bounds, the potential estimate at a fixed $p>2$, and
Sobolev embedding give
$\|\psi_0\|_\infty+\|\theta(\rho_0)_+\|_\infty\le C$,
so $\mathcal M_k(0)\le C_m$.
Gronwall's inequality, \eqref{eq:basic}, and
\eqref{mom:comparison} prove \eqref{mom:moment-bound}.
\end{proof}

Since $k+1\ge4\gamma$, the preceding lemma gives
\begin{equation}\label{mom:pressure-viscosity}
 \sup_{0\le t\le S}
 \left(\|P(t)\|_4+\int_\Om\lambda(\rho(t))\dd\right)\le C_T.
\end{equation}
Since $F=(2\mu+\lambda)\diver u-P$ and $2\mu+\lambda\ge2\mu$,
\eqref{mom:pressure-viscosity} and \eqref{eq:basic} give
\begin{equation}\label{mom:Ecomparison}
 C_T^{-1}(1+A_3^2)\le E\le C_T(1+A_3^2),\qquad
                  \int_0^S E(t)\,\mathrm dt\le C_T.
\end{equation}
Testing momentum by $u$, including its slip boundary term, gives
\[
 A_3^2=\int_\Om P\diver u\dd-\int_\Om\rho\dot u\cdot u\dd
       \le C_T A_3+C_T\sqrt y.
\]
Young's inequality and \eqref{mom:Ecomparison} therefore imply
\begin{equation}\label{mom:Esquare}
                       E^2\le C_T(1+y).
\end{equation}

\subsection{Logarithmic energy estimate}\label{eng:section}

On $[0,S]$, fix
\[
 R=\max\left\{2,\sup_{0\le t\le S}\|\rho(t)\|_\infty\right\}.
\]
Then $h=2\mu+\lambda(\rho)\le CR^b$, where $b\ge1$ is fixed by
\eqref{eq:coeffconsequence}.

\begin{lemma}[Low-exponent elliptic estimate]\label{eng:elliptic-lemma}
For every $t\in[0,S]$,
\begin{equation}\label{eng:ellipticlow}
 \|\nabla F\|_2+\|\nabla\omega\|_2
 \le C_T\left(\sqrt{Ry}+\sqrt E\right).
\end{equation}
Moreover, $|\bar F|\le C_T\sqrt E$ and
$\|F\|_2\le CR^{b/2}\sqrt E$.
\end{lemma}
\begin{proof}
Since $\|\rho\dot u\|_2\le\sqrt{Ry}$ and
$\|u\|_{H^1}\le C_T\sqrt E$ by \eqref{mom:Ecomparison},
\eqref{mom:fluxelliptic} at $p=2$ gives \eqref{eng:ellipticlow}.
Since
$\int_\Om\diver u\dd=0$,
\[
 |\bar F|
 \le C\left|\int_\Om\lambda(\rho)\diver u\dd\right|+C\|P\|_1
 \le C_T\sqrt E.
\]
Here Cauchy--Schwarz uses \eqref{mom:pressure-viscosity}; the bound
$\|F\|_2^2\le CR^bE$ follows from $h\le CR^b$.
\end{proof}

\begin{proposition}[Logarithmic energy]\label{eng:proposition}
There is a positive constant $c$, independent of $R$ and the positive
lower density bound, such that
\begin{equation}\label{eng:differential}
 E'+cy\le C_T R\log(e+R)E^2.
\end{equation}
Consequently,
\begin{equation}\label{eng:logenergy}
 \sup_{0\le t\le S}\log E(t)
 +\int_0^S\frac{e+y}{E}\,\mathrm dt
 \le C_T R\log(e+R).
\end{equation}
\end{proposition}
\begin{proof}
Testing momentum by $2\dot u$ and using
\[
\begin{aligned}
 \diver\dot u&=\Dt(\diver u)+(\diver u)^2-2\det\nabla u,
 &\curl\dot u&=\Dt\omega+\omega\diver u,\\
 \Dt\lambda&=-\rho\lambda'\diver u,
 &\Dt P&=-\gamma P\diver u
\end{aligned}
\]
gives
\begin{equation}\label{eng:identity}
\begin{aligned}
 E'+2y={}&
 \int_\Om\left(1+\frac{\rho\lambda'}h\right)\frac{F^2}h\diver u\dd\\
 &-2\int_\Om\left(1+\frac{\rho\lambda'}h\right)F(\diver u)^2\dd
 +2\gamma\int_\Om\frac{FP}h\diver u\dd\\
 &+4\int_\Om F\det\nabla u\dd
 -\mu\int_\Om\omega^2\diver u\dd\\
 &-2\int_{\partial\Om}\kappa F u_\tau^2\ds
 +\frac{2\mu}{3}\int_{\partial\Om}(\partial_\tau\alpha)u_\tau^3\ds .
\end{aligned}
\end{equation}
The boundary terms follow from $\dot u\cdot n=-\kappa u_\tau^2$,
$\dot u\cdot\tau=\partial_tu_\tau+u_\tau\partial_\tau u_\tau$,
and $\omega=-\alpha u_\tau$:
\[
 2\mu\int_{\partial\Om}\omega(\dot u\cdot\tau)\ds
 =-\mu\frac{\mathrm d}{\mathrm dt}
       \int_{\partial\Om}\alpha u_\tau^2\ds
   +\frac{2\mu}{3}\int_{\partial\Om}(\partial_\tau\alpha)u_\tau^3\ds.
\]
Since $|\rho\lambda'|/h\le C$ and $\diver u=(F+P)/h$, the absolute sum of
the first three terms is at most
\begin{equation}\label{eng:coefficient-terms}
 C\|\diver u\|_2
 \left(\left\|\frac{F^2}h\right\|_2
              +\left\|\frac{FP}h\right\|_2\right).
\end{equation}
For $0<\delta\le1/4$, H\"older's inequality and $h\ge2\mu$ give
\[
 \left\|\frac{F^2}h\right\|_2
 \le C\left(\int_\Om\frac{F^2}h\dd\right)^{(1-\delta)/2}
             \|F\|_{2(1+\delta)/\delta}^{1+\delta}.
\]
For a zero-mean function on $\Om$, the two-dimensional
Gagliardo--Nirenberg inequality has the explicit form
\begin{equation}\label{eng:GN}
 \|f\|_q\le C\sqrt q\,\|f\|_2^{2/q}
                    \|\nabla f\|_2^{1-2/q},\qquad q\ge2.
\end{equation}
Applying it to $F-\bar F$ gives
\[
 \|F-\bar F\|_{2(1+\delta)/\delta}^{1+\delta}
 \le \frac C{\sqrt\delta}
       \|F-\bar F\|_2^\delta\|\nabla F\|_2.
\]
The constant is uniform because
$\delta^{-(1+\delta)/2}\le C\delta^{-1/2}$. Hence
\begin{equation}\label{eng:weighted-quadratic}
 \left\|\frac{F^2}h\right\|_2
 \le \frac{C_T}{\sqrt\delta}R^{b\delta/2}
             \sqrt E\,\|\nabla F\|_2+C_TE.
\end{equation}
For the pressure term, use the fixed bound
$\|P\|_{2+\delta}\le C_T$ and
$q=2(2+\delta)/\delta$ in \eqref{eng:GN}. Weighted arithmetic--geometric
mean then gives
\begin{equation}\label{eng:pressure-quadratic}
 \left\|\frac{FP}h\right\|_2
 \le \frac{C_T}{\sqrt\delta}
       R^{b\delta/[2(2+\delta)]}
       \left(\sqrt E+\|\nabla F\|_2\right).
\end{equation}
Using $\|\diver u\|_2\le C\sqrt E$ and \eqref{eng:ellipticlow} in
\eqref{eng:coefficient-terms}--\eqref{eng:pressure-quadratic}, followed
by Young's inequality, bounds these terms by
\begin{equation}\label{eng:cubic-bound}
 \eps y+\frac{C_{T,\eps}}\delta R^{1+b\delta}E^2 .
\end{equation}

For the Jacobian term, let $p\ge8$. Lemma~\ref{eng:elliptic-lemma}
and \eqref{eng:GN} give
\[
 \|F\|_p+\|\omega\|_p
 \le C_T\sqrt p\,R^{b/p}E^{1/p}
                         (\sqrt E+\sqrt{Ry})^{1-2/p}.
\]
Using \eqref{mom:Hodge}, $\|P\|_p\le CR^\gamma$, and
$\|u\|_2\le C_T\sqrt E$, we obtain
\[
 \|\nabla u\|_p
 \le C_Tp^{13/2}R^{b/p}E^{1/p}(\sqrt E+\sqrt{Ry})^{1-2/p}
                                      +Cp^6R^\gamma.
\]
Interpolating $\|\nabla u\|_2\le C_T\sqrt E$ with this bound at
exponent $2p/(p-1)$, followed by H\"older's inequality, gives
\[
\begin{aligned}
 \left|\int_\Om F\det\nabla u\dd\right|
 &\le C_T\sqrt p
      \left(R^{b/(p-2)}+R^{b/p+2\gamma/(p-2)}\right)\\
 &\qquad\times E(\sqrt E+\sqrt{Ry}).
\end{aligned}
\]
The two powers of $E$ are $1$ and $1-2/[p(p-2)]$, and those of
$\sqrt E+\sqrt{Ry}$ are $1$ and $1-2/p$; both bases are at least one.
Since $p^{13/(p-2)}+p^{12/(p-2)}\le C$ for $p\ge8$, choose
\[
 p=8+(b+2\gamma+1)\log(e+R).
\]
The displayed powers of $R$ are now bounded, so
\begin{equation}\label{eng:Jacobian}
\begin{aligned}
 \left|\int_\Om F\det\nabla u\dd\right|
 &\le C_T\sqrt{\log(e+R)}\,E
                    \left(\sqrt{Ry}+\sqrt E\right)\\
 &\le\eps y+C_{T,\eps}R\log(e+R)E^2 .
\end{aligned}
\end{equation}
For the vorticity term, the bounded-domain inequality
$\|\omega\|_4^2\le
C\|\omega\|_2(\|\nabla\omega\|_2+\|\omega\|_2)$ gives
\[
 \left|\int_\Om\omega^2\diver u\dd\right|
 \le C_TE\left(\sqrt{Ry}+\sqrt E\right)
 \le\eps y+C_{T,\eps}RE^2.
\]
The trace embedding $H^1(\Om)\to L^3(\partial\Om)$ and the controlled
mean of $F$ similarly imply
\[
\begin{aligned}
 \left|\int_{\partial\Om}\kappa F u_\tau^2\ds\right|
 &\le C\|F\|_{H^1}\|u\|_{H^1}^2
 \le C_TE\left(\sqrt{Ry}+\sqrt E\right)
 \le\eps y+C_{T,\eps}RE^2,\\
 \left|\int_{\partial\Om}(\partial_\tau\alpha)u_\tau^3\ds\right|
 &\le C\|u\|_{H^1}^3\le CE^{3/2}.
\end{aligned}
\]
Take $\delta=[4\log(e+R)]^{-1}$ in \eqref{eng:cubic-bound}.
Then $R^{b\delta}\le e^{b/4}$. Substitution of all these estimates into
\eqref{eng:identity} and absorption of $\eps y$ prove
\eqref{eng:differential}.

Dividing that inequality by $E$ gives
\[
 (\log E)'+c\frac yE
 \le C_TR\log(e+R)E.
\]
Since $E(0)\le C$ and $\int_0^S E\,\mathrm dt\le C_T$ by
\eqref{mom:Ecomparison}, integration proves \eqref{eng:logenergy};
also $e/E\le1$.
\end{proof}

\begin{lemma}[High-exponent flux gradient]\label{eng:highgradient-lemma}
For every $m\ge8$,
\begin{equation}\label{eng:highgradient}
 \|\nabla F\|_m
 \le C_Tm^{15/2}R^\gamma
       \left(\sqrt{e+y}+\|\nabla\dot u\|_2\right).
\end{equation}
The constant is independent of $m$ and $R$.
\end{lemma}
\begin{proof}
The mass-weighted Poincar\'e inequality and the two-dimensional
Sobolev inequality give
\[
 \|\dot u\|_m
 \le C_T\sqrt m\left(\sqrt y+\|\nabla\dot u\|_2\right),
 \qquad
 \|\rho\dot u\|_m
 \le C_TR\sqrt m\left(\sqrt y+\|\nabla\dot u\|_2\right).
\]
The Poincar\'e constant depends only on the mass and the basic
$L^\gamma$ density bound. By \eqref{eng:ellipticlow} and Sobolev,
\[
 \|F\|_m+\|\omega\|_m
 \le C_T\sqrt m\left(\sqrt{Ry}+\sqrt E\right).
\]
Using \eqref{mom:Hodge}, $h\ge2\mu$, and
$\|P\|_m\le CR^\gamma$, we obtain
\[
\begin{aligned}
 \|\nabla u\|_m
 &\le C_Tm^{13/2}\left(\sqrt{Ry}+\sqrt E\right)+Cm^6R^\gamma
                          +C\|u\|_2\\
 &\le C_Tm^{13/2}R^\gamma\sqrt{e+y},
 \qquad
 \|u\|_m\le C_T\sqrt m\,\sqrt{e+y}.
\end{aligned}
\]
The last line uses \eqref{mom:Esquare} and $\gamma>1$.
Finally the quantitative elliptic estimate \eqref{mom:fluxelliptic} gives
\[
 \|\nabla F\|_m
 \le Cm\bigl(\|\rho\dot u\|_m
                     +\|\nabla u\|_m+\|u\|_m\bigr).
\]
This proves \eqref{eng:highgradient}.
\end{proof}

\subsection{Acceleration estimate}\label{acc:section}

We correct the normal acceleration to obtain a tangent test field.

\begin{lemma}[Positive-time acceleration]\label{acc:lemma}
There is a fixed $N\ge1$, depending only on $\gamma$, the coefficient
and domain constants, such that, for $0<\tau<S\le T$,
\begin{equation}\label{acc:logacceleration}
 \sup_{\tau\le t\le S}\log(e+A_2^2(t))
 +\int_\tau^S
 \left(1+\frac{\|\nabla\dot u\|_2}{\sqrt{e+A_2^2}}\right)^2\,\mathrm dt
 \le C_{\tau,T}R^N.
\end{equation}
The constant is independent of the positive density lower bound and of
$A_2(0)$.
\end{lemma}

\begin{proof}
Keep $R$ fixed on $[0,S]$. The fixed exponent
$N$ may be enlarged in the estimates below.

By \eqref{mom:Hodge}, \eqref{eng:GN}, \eqref{eng:ellipticlow}, and
\eqref{mom:pressure-viscosity},
\begin{equation}\label{acc:gradientfour}
 \|\nabla u\|_4^4
 \le C_T R^N E(y+E)\le C_T R^N E(1+y),
\end{equation}
where the last inequality follows from \eqref{mom:Esquare}.

Extend $n$ constantly along normal segments in a boundary collar and
fix a bounded $C^{1,1}$ extension to $\overline\Om$. Set
\[
 q=\bigl((u\cdot\nabla)n\cdot u\bigr)n,\qquad
 v=\dot u+q,\qquad z=\int_\Om\rho|v|^2\dd.
\]
On the boundary,
\[
 q=\kappa u_\tau^2n,\qquad
 \dot u\cdot n=-\kappa u_\tau^2,\qquad
 v\cdot n=0,\qquad v_\tau=\dot u\cdot\tau.
\]
The extension is fixed, so
\begin{equation}\label{acc:qbounds}
\begin{aligned}
 |q|&\le C|u|^2,&
 |\nabla q|&\le C\bigl(|u|^2+|u||\nabla u|\bigr),&
 |\Dt q|&\le C\bigl(|u||\dot u|+|u|^3\bigr),\\
 \int_\Om\rho|q|^2\dd&\le C_T R E^2,&
 y&\le2z+C_T R E^2,&
 z&\le2y+C_T R E^2.
\end{aligned}
\end{equation}
By \eqref{acc:gradientfour}, \eqref{mom:Esquare}, and Young's inequality,
\begin{equation}\label{acc:qgradient}
\begin{aligned}
 \|\nabla q\|_2^2
 &\le C_T\bigl(E^2+E\|\nabla u\|_4^2\bigr)\\
 &\le C_T R^N\bigl(E^2+E^{3/2}(1+y)^{1/2}\bigr)
 \le C_T R^N E(1+y).
\end{aligned}
\end{equation}
To differentiate momentum, use $(\nabla u)_{ij}=\partial_j u_i$ and set
\[
 \mathsf S=F I+\mu\omega J,\qquad
 \Sigma=\Dt\mathsf S+(\diver u)\mathsf S
                         -\mathsf S(\nabla u)^{\mathsf T}.
\]
Since $\diver\mathsf S=\rho\dot u$, the continuity equation gives
\begin{equation}\label{acc:commutation}
 \rho\Dt\dot u_i=\partial_j\Sigma_{ij}.
\end{equation}
Since $F=h\diver u-P$,
\[
\begin{aligned}
 \Dt F&=h\diver\dot u-\rho\lambda'(\diver u)^2
          -h\operatorname{tr}((\nabla u)^2)+\gamma P\diver u,\\
 \Dt\omega&=\curl\dot u-\omega\diver u.
\end{aligned}
\]
Consequently,
\begin{equation}\label{acc:stressremainder}
\begin{aligned}
 \Sigma&=h(\diver\dot u)I+\mu(\curl\dot u)J+Q,\\
 Q&=\bigl((h-\rho\lambda')(\diver u)^2
          -h\operatorname{tr}((\nabla u)^2)
          +(\gamma-1)P\diver u\bigr)I\\
  &\qquad-F(\nabla u)^{\mathsf T}-\mu\omega J(\nabla u)^{\mathsf T},\\
 |Q|&\le C\bigl(h|\nabla u|^2+P|\nabla u|\bigr).
\end{aligned}
\end{equation}
Choose a fixed sufficiently large $K$ and define
\[
 H=1+z+K R E^2.
\]
Equations~\eqref{mom:Esquare} and \eqref{acc:qbounds} imply
\begin{equation}\label{acc:Hcomparison}
 1+y\le C_T H,\qquad H\le C_T R(1+y).
\end{equation}
Test $\rho\Dt v=\diver\Sigma+\rho\Dt q$ by $v$ and set
$\mathcal D_v=\int_\Om(h(\diver v)^2+\mu(\curl v)^2)\dd$.
Since $v\cdot n=0$, \eqref{ell:weightedHodge} gives
\[
 \|v\|_{H^1}^2\le C(\mathcal D_v+z),\qquad
 \mathcal D_v\ge c\|\nabla v\|_2^2-Cz.
\]
Thus $\eps\|\nabla v\|_2^2\le C\eps\mathcal D_v+C\eps z$.
The lower order terms satisfy
\[
 \sqrt{Ry}\sqrt E\sqrt z\le CR E(y+z),\qquad
 \sqrt z(E+\sqrt E\|\nabla u\|_4)
 \le z+C(E^2+E\|\nabla u\|_4^2).
\]
These are bounded by $CR^NEH$ since $z\le H$ and $E\ge e$.
The cross terms containing $q$ are bounded by a small multiple of the
principal dissipation plus $C R^b\|\nabla q\|_2^2$. Also,
\[
\begin{aligned}
 \left|\int_\Om Q:\nabla v\dd\right|
 &\le\eps\|\nabla v\|_2^2
   +C_\eps R^{2b}\|\nabla u\|_4^4
   +C_\eps R^{2\gamma}\|\nabla u\|_2^2\\
 &\le\eps\mathcal D_v+C_{\eps,T}R^N E H.
\end{aligned}
\]
The correction in the material derivative satisfies
\[
\begin{aligned}
 \int_\Om\rho|u||\dot u||v|\dd
 &\le\sqrt{Ry}\,\|u\|_4\|v\|_4
 \le\eps\mathcal D_v+C_{\eps,T}R E(y+z),\\
 \int_\Om\rho|u|^3|v|\dd
 &\le\sqrt{Rz}\,\|u\|_6^3
 \le C_T\sqrt{Rz}\,E^{3/2}
 \le C_T E(z+R E^2).
\end{aligned}
\]

For the boundary traction, use $\omega=-\alpha u_\tau$ and
$\Dt u_\tau=v_\tau$ on $\partial\Om$:
\[
\begin{aligned}
 (\nabla u)^{\mathsf T}n
  &=(\partial_nu\cdot n)n-\kappa u_\tau\tau,\\
 \diver u&=\partial_nu\cdot n+\partial_\tau u_\tau,\\
 \Dt\omega&=-\alpha v_\tau-(\partial_\tau\alpha)u_\tau^2.
\end{aligned}
\]
Using $Jn=\tau$ and $v\cdot n=0$ gives
\begin{equation}\label{acc:boundarytraction}
\begin{aligned}
 (\Sigma n)\cdot v
 &=\mu\bigl(\Dt\omega+\omega\partial_\tau u_\tau\bigr)v_\tau
                                      +\kappa F u_\tau v_\tau\\
 &=-\mu\alpha v_\tau^2
   -\mu(\partial_\tau\alpha)u_\tau^2v_\tau
   -\mu\alpha u_\tau(\partial_\tau u_\tau)v_\tau
   +\kappa F u_\tau v_\tau.
\end{aligned}
\end{equation}
The first term is nonpositive, and the trace inequalities give
\[
\begin{aligned}
 \left|\int_{\partial\Om}\kappa F u_\tau v_\tau\ds\right|
 &\le C\|F\|_{H^1}\|u\|_{H^1}\|v\|_{H^1}\\
 &\le\eps\mathcal D_v+C_{\eps,T}R^N E H,\\
 \left|\int_{\partial\Om}(\partial_\tau\alpha)u_\tau^2v_\tau\ds\right|
 &\le C_T E(\|\nabla v\|_2+\sqrt z)
 \le\eps\mathcal D_v+C_{\eps,T}R^N E H.
\end{aligned}
\]
For the remaining term, use fixed smooth extensions of $\alpha,\tau$
and write $u_\tau=u\cdot\tau$, $v_\tau=v\cdot\tau$ in $\Om$.
Stokes' formula gives
\[
\begin{aligned}
 \left|\int_{\partial\Om}\alpha u_\tau(\partial_\tau u_\tau)v_\tau\ds\right|
 &=\frac12\left|\int_{\partial\Om}\alpha v_\tau
                                      \partial_\tau(u_\tau^2)\ds\right|\\
 &\le C\int_\Om|\nabla(\alpha v_\tau)|\,|\nabla(u_\tau^2)|\dd\\
 &\le C_T(\|\nabla v\|_2+\sqrt z)
                  \bigl(E+\sqrt E\,\|\nabla u\|_4\bigr)\\
 &\le\eps\mathcal D_v
       +C_{\eps,T}\bigl(z+E^2+E\|\nabla u\|_4^2\bigr)\\
 &\le\eps\mathcal D_v+C_{\eps,T}R^N E H.
\end{aligned}
\]
Combining these estimates and absorbing $\eps\mathcal D_v$ gives
\begin{equation}\label{acc:zdifferential}
 z'+c\|\nabla v\|_2^2\le C_T R^N E H.
\end{equation}
By \eqref{eng:differential},
\[
 E'+cy\le C_T R^N E^2.
\]
Multiply by $2KRE$ and add \eqref{acc:zdifferential}.
Since $KRE^3\le EH$, dropping the nonnegative $Ey$ term gives
\begin{equation}\label{acc:Hdifferential}
 H'+c\|\nabla v\|_2^2\le C_T R^N E H.
\end{equation}
Here $\int_0^S E\,\mathrm dt\le C_T$ by \eqref{mom:Ecomparison}.

Finally, \eqref{eng:logenergy} implies, after enlarging $N$ if needed,
\[
 \int_0^S y\,\mathrm dt\le C_T\exp(C_T R^N).
\]
There exists $s_\tau\in[\tau/2,\tau]$ such that
\[
 y(s_\tau)\le\frac2\tau\int_0^\tau y\,\mathrm dt,
 \qquad \log H(s_\tau)\le C_{\tau,T}R^N,
\]
where the second bound uses \eqref{acc:Hcomparison}.
Divide \eqref{acc:Hdifferential} by $H$ and integrate from $s_\tau$.
Using again \eqref{mom:Ecomparison}, we find
\[
 \sup_{\tau\le t\le S}\log H(t)
       +\int_\tau^S\frac{\|\nabla v\|_2^2}{H}\,\mathrm dt
 \le C_{\tau,T}R^N.
\]
Since $\dot u=v-q$, \eqref{acc:qgradient} and
\eqref{acc:Hcomparison} give
\[
 \frac{\|\nabla\dot u\|_2^2}{e+y}
 \le C_T R\frac{\|\nabla v\|_2^2}{H}+C_T R^N E.
\]
Integration, the first inequality in \eqref{acc:Hcomparison}, and
$(1+x)^2\le2+2x^2$ prove \eqref{acc:logacceleration}.
\end{proof}

\subsection{Density upper bound}

For $\bar f=0$ and $p,m\ge8$,
\begin{equation}\label{cl:interpolation}
 \|f\|_\infty\le C\sqrt p\,
 \|f\|_2^\eta\|\nabla f\|_2^{1-\eta-\sigma}
                         \|\nabla f\|_m^\sigma,\quad
 \eta=\frac{2(m-2)}{p(m-2)+2m},\quad
 \sigma=\frac{2m}{p(m-2)+2m}.
\end{equation}
The constant is independent of $p,m$. Take a fixed reflection extension
with $\|\mathcal Ef\|_q\le C\|f\|_q$ and
$\|\nabla\mathcal Ef\|_q\le C\|\nabla f\|_q$, $q\ge2$.
Here $C$ is uniform in $q$, and Poincare's inequality absorbs the cutoff
terms. Balance
$\ell^{-2/p}\|\mathcal Ef\|_p$ and
$\ell^{1-2/m}\|\nabla\mathcal Ef\|_m$ and apply \eqref{eng:GN}.
Poincare--Morrey handles radii above the domain scale. This proves
\eqref{cl:interpolation}.

\begin{lemma}[Weighted and local flux estimates]\label{cl:flux}
For $0<\tau<S\le T$,
\begin{equation}\label{cl:weightedflux}
 \int_\tau^S\frac{\|F\|_\infty^2}{E}\,\mathrm dt
       \le C_{\tau,T}R_S^2[\log(e+R_S)]^2.
\end{equation}
For $\tau\le s<t\le S$,
\begin{equation}\label{cl:localflux}
 \int_s^t\|F\|_\infty\,\mathrm d\xi
 \le C_{\tau,T}R_t\log(e+R_t)
       \left(\int_s^t(1+A_3^2)\,\mathrm d\xi\right)^{1/2}.
\end{equation}
\end{lemma}
\begin{proof}
Freeze $R=R_S$ and write
\[
 w=1+\frac{e+y}{E},\qquad
 Z=1+\frac{\|\nabla\dot u\|_2}{\sqrt{e+y}}.
\]
From \eqref{eng:ellipticlow},
$\|F\|_2\le CR^{b/2}\sqrt E$,
\eqref{eng:highgradient}, and \eqref{cl:interpolation},
\[
 \frac{\|F-\bar F\|_\infty}{\sqrt E}
 \le C_T\sqrt p\,
 R^{1/2+(b-1)\eta/2+(\gamma-1/2)\sigma}
          m^{15\sigma/2}w^{(1-\eta)/2}Z^\sigma.
\]
Choose
\[
 p=8+(N+b+2\gamma+2)\log(e+R),\qquad m=e^p,
\]
where $N$ is enlarged to cover \eqref{acc:logacceleration}.
Then
\[
 p\eta\le2,\quad p\sigma\le\frac83,\quad
 \sigma-\eta=\frac4{p(m-2)+2m}.
\]
The extra powers of $R$ and $m^{15\sigma/2}$ are bounded. Moreover,
\[
 \sup_{[\tau,S]}w^{(\sigma-\eta)/2}
 \le\exp\left(\frac{C_{\tau,T}R^N}{pm}\right)\le C_{\tau,T},
\]
by \eqref{acc:logacceleration} and $E\ge e$.
Consequently
\[
 \frac{\|F-\bar F\|_\infty}{\sqrt E}
       \le C_{\tau,T}\sqrt{pR}\,
                 w^{(1-\sigma)/2}Z^\sigma.
\]
Holder in time, \eqref{eng:logenergy}, and
\eqref{acc:logacceleration} give
\[
 \left(\int_\tau^S\frac{\|F-\bar F\|_\infty^2}{E}\,\mathrm dt\right)^{1/2}
 \le C_{\tau,T}\sqrt p\,
 R^{1+(N-1)\sigma/2}[\log(e+R)]^{(1-\sigma)/2}
 \le C_{\tau,T}R\log(e+R).
\]
Since $|\bar F|\le C_T\sqrt E$ by Lemma~\ref{eng:elliptic-lemma},
\eqref{cl:weightedflux} follows. Apply it on $[0,t]$; Cauchy--Schwarz
on $[s,t]$ and \eqref{mom:Ecomparison} give \eqref{cl:localflux}.
\end{proof}

\begin{lemma}[A uniform initial interval]\label{cl:initial}
There are $\tau_0>0$ and $M_0$, depending only on the data in
Proposition~\ref{thm:main}, such that
\[
 \sup_{0\le t\le\min\{S,\tau_0\}}\|\rho(t)\|_\infty\le M_0.
\]
\end{lemma}
\begin{proof}
Work before first contact with $R_*=2(1+C_0)$, with $t\le1$.
By \eqref{eng:logenergy} and \eqref{acc:Hcomparison},
\[
 \sup_{s\le t}E(s)+\int_0^t(y+H)\,\mathrm ds\le C_*.
\]
Multiplying \eqref{acc:Hdifferential} by time and using
\eqref{acc:qgradient} gives
\[
 \sup_{0<s\le t}sy(s)
       +\int_0^t s\|\nabla\dot u\|_2^2\,\mathrm ds\le C_*.
\]
The weighted initial term vanishes for each smooth solution.
By \eqref{eng:highgradient}
with $m=8$,
$\|\nabla F\|_4\le C_*(\sqrt{e+y}+\|\nabla\dot u\|_2)$.
Since $\|F\|_2\le C_*$, the GN inequality gives
\[
 \|F\|_\infty
 \le C_*\bigl(1+y^{1/3}+\|\nabla\dot u\|_2^{2/3}\bigr).
\]
H\"older's inequality therefore yields
\[
\begin{aligned}
 \int_0^t\|F\|_\infty\,\mathrm ds
 &\le C_*\left[t+t^{2/3}\left(\int_0^t y\,\mathrm ds\right)^{1/3}
       +t^{1/3}\left(\int_0^t
                   s\|\nabla\dot u\|_2^2\,\mathrm ds\right)^{1/3}\right]\\
 &\le C_*t^{1/3}.
\end{aligned}
\]
Since $u\cdot n=0$, the flow preserves $\overline\Om$, and
\begin{equation}\label{cl:thetaflow}
 \Dt\theta(\rho)+P=-F.
\end{equation}
Thus $\theta(\rho(t))\le\theta(C_0)+C_*t^{1/3}$.
Choose $\tau_0$ with
$\theta(C_0)+C_*\tau_0^{1/3}<\theta(R_*)$ to exclude first contact.
\end{proof}

\begin{proof}[Proof of Proposition~\ref{thm:main}]
If $S\le\tau_0$, use Lemma~\ref{cl:initial}. Otherwise fix
$\tau=\tau_0/2$ and
$M_*>\max\{2,r_\lambda,M_0\}$ with $\theta(M_*)>0$.
Define
\[
 r_j=\theta^{-1}(2^j\theta(M_*)),\qquad t_0=\tau,
\]
and let $t_j$, $j\ge1$, be the first time the running maximum reaches
$r_j$, considering only attained levels. Write
\[
 \mathcal E_j=\int_{t_{j-1}}^{t_j}(1+A_3^2)\,\mathrm dt.
\]
For a trajectory ending at density $r_j$, integration of
\eqref{cl:thetaflow} and \eqref{cl:localflux} give
\[
 \frac{\theta(r_j)}2
 \le C_{\tau,T}r_j\log(e+r_j)\sqrt{\mathcal E_j}.
\]
Since $\theta(r_j)\ge cr_j[\log(e+r_j)]^a$,
$r_j\le C2^j\theta(M_*)$, and $a\le1$, every finite number $L$
of attained levels satisfies
\[
\begin{aligned}
 C_T\ge\sum_{j=1}^{L}\mathcal E_j
 &\ge c_T\sum_{j=1}^{L}
       \left[\frac{\theta(r_j)}{r_j\log(e+r_j)}\right]^2\\
 &\ge c_T\sum_{j=1}^{L}(j+1)^{2a-2}.
\end{aligned}
\]
The last sum diverges for $a>1/2$, so $R_S\le C_T$.

Now \eqref{eng:logenergy} gives $\sup E+\int_0^S y\,\mathrm dt\le C_T$,
and \eqref{mom:Ecomparison} controls $A_3$ and $u$.
Equation \eqref{acc:Hdifferential} has a bounded coefficient and,
by \eqref{acc:Hcomparison}, $\int_0^S H\,\mathrm dt\le C_T$.
Multiplying by time and integrating yields
$\sup tH+\int_0^S t\|\nabla v\|_2^2\,\mathrm dt\le C_T$.
Equations \eqref{acc:Hcomparison} and \eqref{acc:qgradient}
give the remaining terms of \eqref{eq:main}.
\end{proof}

\subsection{Uniform density moments and pressure dissipation}

Assume that $\alpha\ge0$ and $\alpha\not\equiv0$. Set
\[
 r_*=\bar\rho=\frac{M}{|\Om|},
\]
and
\[
 A_1^2=\int_\Om\left(\mu\omega^2+\frac{G^2}{h}\right)\dd
                  +\mu\int_{\partial\Om}\alpha u_\tau^2\ds.
\]
All constants in this subsection are independent of $S>0$ and the
positive density lower bound. By \eqref{uni:coercive} and
\eqref{eq:basic},
\begin{equation}\label{uni:basic}
 \|u\|_{H^1}\le CA_3,\qquad
 \sup_{t\le S}\int_\Om(\rho|u|^2+\rho^\gamma)\dd
                    +\int_0^S A_3^2\,\mathrm dt\le C.
\end{equation}

\begin{lemma}\label{uni:moment-lemma}
For every integer $k\ge3$ with $(2a-1)k\ge1$,
\begin{equation}\label{uni:moment}
 \sup_{t\le S}\int_\Om\left(
       \rho\theta(\rho)_+^k+
       \rho^{k+1}[\log(e+\rho)]^{ak}\right)\dd\le C_k.
\end{equation}
Consequently, $\sup_{t\le S}\|\rho(t)\|_p\le C_p$ for every fixed
$1\le p<\infty$, and
\begin{equation}\label{uni:momentum}
 \|\rho u\|_m\le C_m A_3,\qquad 2\le m<\infty.
\end{equation}
For every fixed $p\ge2$, there are $L_p\ge2(1+r_*)$ and $C_p$ such that
\begin{equation}\label{uni:tails}
 \int_0^S\int_\Om(\rho-L_p)_+^p\dd\,\mathrm dt\le C_p.
\end{equation}
\end{lemma}

\begin{proof}
Let $\psi=Q(\rho u)$ and $B=Q(\nabla^\perp\omega)$. Since
$G=F+\bar P$ and $\bar G=|\Om|^{-1}\int_\Om\lambda\diver u\dd$,
\eqref{mom:transport} becomes
\begin{equation}\label{uni:transport}
 \Dt(\theta(\rho)-\psi)+P
       =\bar P-\bar G-\mu B-\mathfrak C(u,\rho u).
\end{equation}
Fix a shift $L>0$, to be chosen below, and put
\[
 \varphi=(\theta(\rho)-\psi-L)_+,\qquad
 \mathcal M_k=\int_\Om\rho\varphi^k\dd,\qquad m=k+1.
\]
The potential estimate and the kinetic energy bound give
\[
 \|\psi\|_{2m}\le C_k\|\rho u\|_{2m/(m+1)}
                    \le C_k\|\rho\|_m^{1/2},\qquad
 \int_\Om\rho|\psi|^k\dd\le C_k\|\rho\|_m^{m/2}.
\]
Using $\theta_+\le\varphi+|\psi|+L$, the growth condition, and
absorption of the resulting square root, we obtain
\begin{equation}\label{uni:moment-comparison}
\begin{aligned}
 &\int_\Om\left(\rho\theta(\rho)_+^k+
              \rho^m[\log(e+\rho)]^{ak}\right)\dd
                      \le C_{k,L}(1+\mathcal M_k),\\
 &\int_\Om\rho|\psi|^k\dd
                      \le C_{k,L}(1+\mathcal M_k)^{1/2}.
\end{aligned}
\end{equation}

By \eqref{uni:basic} and Trudinger's inequality,
\[
 \int_\Om\exp\left(\frac{c|u|^2}{A_3^2}\right)\dd\le C
 \quad\text{if }A_3>0,\qquad u=0\quad\text{if }A_3=0.
\]
Orlicz--Young and $ak\ge m/2$ therefore give
\[
 \|\rho u\|_m\le C_{k,L}A_3(1+\mathcal M_k)^{1/m}.
\]
Together with Lemmas~\ref{mom:commutator} and \ref{mom:boundary}, this
implies
\[
 \|\mathfrak C(u,\rho u)\|_m
       \le C_{k,L}A_3^2(1+\mathcal M_k)^{1/m},\qquad
 \|B\|_{2m}\le C_kA_3.
\]
For any $Z$, H\"older's inequality gives
\[
 \int_\Om\rho\varphi^{k-1}|Z|\dd
 \le\mathcal M_k^{1-1/k}\|\rho\|_m^{1/k}\|Z\|_m.
\]
Since $1-1/k+1/(km)+1/m=1$, the commutator contribution is bounded
by $C_{k,L}A_3^2(1+\mathcal M_k)$. Moreover,
\[
\begin{aligned}
 |\bar G|^2
 &\le CA_3^2\int_\Om\lambda(\rho)\dd
 \le C_{k,L}A_3^2(1+\mathcal M_k)^{1/k},\\
 \mu\int_\Om\rho\varphi^{k-1}|B|\dd
 &\le C_{k,L}A_3^2(1+\mathcal M_k)
                         +C\int_\Om\rho\varphi^{k-1}\dd.
\end{aligned}
\]
Here the second estimate uses $\mu|B|\le\mu^2 B^2+1$ and
$\|B^2\|_m\le C_kA_3^2$. Similarly,
$|\bar G|\le|\bar G|^2+1$. Testing \eqref{uni:transport} thus gives
\[
 \mathcal M_k'+k\int_\Om\rho^{\gamma+1}\varphi^{k-1}\dd
 \le C_{k,L}A_3^2(1+\mathcal M_k)
                         +Ck\int_\Om\rho\varphi^{k-1}\dd,
\]
where $C$ is independent of $L$; the bound for $\bar P$ follows from
\eqref{uni:basic}. Choose
\[
 L>1+\max\{\theta((4C)^{1/\gamma}),0\}.
\]
On $\{\rho^\gamma\ge4C\}$, half the pressure term absorbs the last
integral. On its complement, $\varphi\le|\psi|$, and
\[
 \|\psi\|_\infty
 \le C\|\rho u\|_{8/3}
 \le C\|\rho\|_4\|u\|_8
 \le C_{k,L}A_3(1+\mathcal M_k)^{1/(k+1)}.
\]
Since $k\ge3$, \eqref{uni:moment-comparison} gives
\[
\begin{aligned}
 \int_{\{\rho^\gamma<4C\}}\rho\varphi^{k-1}\dd
 &\le\|\psi\|_\infty^2\int_\Om\rho|\psi|^{k-3}\dd\\
 &\le C_{k,L}A_3^2(1+\mathcal M_k)^{2/(k+1)+1/2}
 \le C_{k,L}A_3^2(1+\mathcal M_k).
\end{aligned}
\]
Consequently,
\begin{equation}\label{uni:moment-differential}
 \mathcal M_k'+\frac k2
       \int_\Om\rho^{\gamma+1}\varphi^{k-1}\dd
       \le C_{k,L}A_3^2(1+\mathcal M_k).
\end{equation}
The initial bound $\mathcal M_k(0)\le C_k$ and
\eqref{uni:basic} imply
\begin{equation}\label{uni:moment-damping}
 \sup_{t\le S}\mathcal M_k(t)
       +\int_0^S\int_\Om\rho^{\gamma+1}\varphi^{k-1}\dd\,\mathrm dt
       \le C_k.
\end{equation}
This proves \eqref{uni:moment}. Taking $k$ arbitrarily large gives all
fixed density moments, and H\"older's inequality with
\eqref{uni:basic} proves \eqref{uni:momentum}.

For $p\ge3$, choose $k\ge p$ above and increase $L_p$ so that
$\theta(r)\ge2L$ and $\theta(r)\ge cr$ for $r\ge L_p$.
On $\{\rho\ge L_p\}$, either $\varphi\ge\theta(\rho)/4$ or
$\psi\ge\theta(\rho)/4\ge c\rho$. Hence
\[
\begin{aligned}
 \int_{\{\rho\ge L_p\}}\rho^p\dd
 &\le C_p\int_\Om\rho^{\gamma+1}\varphi^{k-1}\dd
       +C_p\|\psi\|_\infty^2\int_\Om\rho|\psi|^{p-3}\dd\\
 &\le C_p\int_\Om\rho^{\gamma+1}\varphi^{k-1}\dd+C_pA_3^2.
\end{aligned}
\]
Integration and \eqref{uni:moment-damping} give \eqref{uni:tails}.
For $2\le p<3$, increase the threshold by one and use
$(\rho-L_3-1)_+^p\le(\rho-L_3)_+^3$.
\end{proof}

\begin{lemma}\label{uni:pressure-lemma}
The relative pressure satisfies
\begin{equation}\label{uni:pressure}
\begin{aligned}
 &\int_0^S\int_\Om(\rho+r_*)^{\gamma-1}(\rho-r_*)^2\dd\,\mathrm dt\\
 &\qquad+\int_0^S\bigl(\|P-r_*^\gamma\|_2^2
                         +\|P-r_*^\gamma\|_4^4\bigr)\,\mathrm dt
 \le C.
\end{aligned}
\end{equation}
\end{lemma}

\begin{proof}
Let
\[
 \Delta\phi=\rho-r_*,\qquad \partial_n\phi=0,\qquad\bar\phi=0.
\]
Then $\phi_t=Q(\rho u)$. Testing conservative momentum by
$\nabla\phi$, which is tangent, gives
\begin{equation}\label{uni:pressure-test}
\begin{aligned}
 \int_\Om(P-r_*^\gamma)(\rho-r_*)\dd
 ={}&\frac{\mathrm d}{\mathrm dt}
          \int_\Om\rho u\cdot\nabla\phi\dd
       +\|\nabla Q(\rho u)\|_2^2\\
 &-\int_\Om\rho u\otimes u:D^2\phi\dd
       +\int_\Om h\diver u\,(\rho-r_*)\dd
       +\mu\int_{\partial\Om}\alpha u_\tau\partial_\tau\phi\ds.
\end{aligned}
\end{equation}
The last term is the Navier boundary contribution; the vorticity part
vanishes since $\curl\nabla\phi=0$. The uniform moments, fixed-exponent
elliptic estimates, and \eqref{uni:basic} give
\[
\begin{aligned}
 \|\nabla Q(\rho u)\|_2^2
       +\left|\int_\Om\rho u\otimes u:D^2\phi\dd\right|
 &\le CA_3^2,\\
 \left|\int_\Om\rho u\cdot\nabla\phi\dd\right|
 &\le\|\sqrt\rho u\|_2\|\rho\|_2^{1/2}\|\nabla\phi\|_4\le C,\\
 \left|\mu\int_{\partial\Om}\alpha u_\tau\partial_\tau\phi\ds\right|
 &\le CA_3\|\rho-r_*\|_2
 \le\eps\|\rho-r_*\|_2^2+C_\eps A_3^2.
\end{aligned}
\]
Split the viscous cross term at $\rho=2L_{b+2}$. By
\eqref{eq:coeffconsequence},
\[
 \left|\int_\Om h\diver u\,(\rho-r_*)\dd\right|
 \le\eps\|\rho-r_*\|_2^2+C_\eps A_3^2
              +C\int_\Om(\rho-L_{b+2})_+^{b+2}\dd.
\]
Since
\[
 (\rho^\gamma-r_*^\gamma)(\rho-r_*)
       \ge c(\rho+r_*)^{\gamma-1}(\rho-r_*)^2,
\]
integration of \eqref{uni:pressure-test}, with $\eps$ small, and
\eqref{uni:tails} give the first bound in \eqref{uni:pressure}.
For $p=2,4$, splitting at $\rho=2L_{p\gamma}$ gives
\[
 |\rho^\gamma-r_*^\gamma|^p
       \le C_p(\rho-r_*)^2+C_p(\rho-L_{p\gamma})_+^{p\gamma}.
\]
The first bound and \eqref{uni:tails} prove the remaining terms.
\end{proof}

Expanding $G=h\diver u-(P-\bar P)$ and using the uniform pressure
moment gives
\begin{equation}\label{uni:comparison}
\begin{aligned}
 A_1^2&\le C\bigl(A_3^2+\|P-\bar P\|_2^2\bigr),\qquad
 A_3^2\le C\bigl(A_1^2+\|P-\bar P\|_2^2\bigr),\\
 C^{-1}(1+A_3^2)&\le e+A_1^2\le C(1+A_3^2),\qquad
             \int_0^S A_1^2\,\mathrm dt\le C.
\end{aligned}
\end{equation}
In particular, since $\|P-\bar P\|_p\le2\|P-r_*^\gamma\|_p$,
\begin{equation}\label{uni:dissipation}
 D(t)=A_3^2+\|P-\bar P\|_2^2+\|P-\bar P\|_4^4,\qquad
                     \int_0^S D(t)\,\mathrm dt\le C.
\end{equation}

\begin{lemma}\label{uni:conditional-decay-lemma}
Fix $R=R_S$. There are $C,c>0$ and a fixed $N\ge1$, independent of
$S$ and $R$, such that
\begin{equation}\label{uni:conditionaldecay}
 \|\rho(t)-r_*\|_2^2
       \le CR^N\exp(-cR^{-N}t),\qquad 0\le t\le S.
\end{equation}
Moreover, uniformly for $m\ge8$,
\begin{equation}\label{uni:pressure-high}
 \int_0^S\|P-\bar P\|_m^2\,\mathrm dt\le CmR^N.
\end{equation}
\end{lemma}

\begin{proof}
With $\phi$ as in \eqref{uni:pressure-test}, set
\begin{equation}\label{uni:relative-energy}
\begin{aligned}
 \mathcal H&=\int_\Om\left[
       \frac12\rho|u|^2+
       \frac{\rho^\gamma-r_*^\gamma
                    -\gamma r_*^{\gamma-1}(\rho-r_*)}{\gamma-1}
                         \right]\dd,\\
 K&=\int_\Om\rho u\cdot\nabla\phi\dd.
\end{aligned}
\end{equation}
Then $\mathcal H'+A_3^2=0$. Choose a fixed $N_0>b+\gamma+2$.
The relative internal energy, $0\le\rho\le R$, and
\eqref{uni:coercive} imply
\[
 cR^{-N_0}\|\rho-r_*\|_2^2\le\mathcal H
       \le CR^{N_0}(A_3^2+\|\rho-r_*\|_2^2),\qquad
 |K|\le CR^{N_0}\mathcal H.
\]
For the last bound, use
$|K|\le C\|\sqrt\rho u\|_2\|\rho-r_*\|_2$, which follows from
the uniform $L^2$ density moment and the Neumann $H^2$ estimate.
In \eqref{uni:pressure-test}, the viscous term now satisfies
\[
 \left|\int_\Om h\diver u\,(\rho-r_*)\dd\right|
       \le\eps\|\rho-r_*\|_2^2+C_\eps R^bA_3^2.
\]
The projection, convection, and boundary estimates above therefore
give
\[
 K'\ge c\|\rho-r_*\|_2^2-CR^bA_3^2.
\]
Take $\delta=c_0R^{-N_0}$, with $c_0>0$ sufficiently small. Then
\[
 \mathcal H-\delta K\asymp\mathcal H,\qquad
 (\mathcal H-\delta K)'
       +c(A_3^2+\delta\|\rho-r_*\|_2^2)\le0,
\]
and hence
\[
 (\mathcal H-\delta K)'
       +cR^{-2N_0}(\mathcal H-\delta K)\le0.
\]
Since $\mathcal H(0)\le C$, this proves
\eqref{uni:conditionaldecay} with a fixed enlarged exponent $N$.

For $m\ge8$, interpolation and $\|\rho-r_*\|_\infty\le CR$ give
\[
 \|P-\bar P\|_m^2
       \le CR^{2\gamma}\|\rho-r_*\|_2^{4/m}
       \le CR^{2\gamma+N/4}\exp(-2cR^{-N}t/m).
\]
Integration gives \eqref{uni:pressure-high}, after enlarging $N$ once
more. All constants and exponents are independent of $m$ and $S$.
\end{proof}

\subsection{Uniform logarithmic energy and acceleration}

Fix $S>0$ and
\[
 R=\max\{2,\sup_{0\le t\le S}\|\rho(t)\|_\infty\},\qquad
 E_*=e+A_1^2.
\]
Constants in this subsection are independent of $S$, $R$, and the
positive density lower bound. We use $D(t)$ from \eqref{uni:dissipation}.

\begin{lemma}[Uniform logarithmic energy]\label{uni:logenergy-lemma}
One has
\begin{equation}\label{uni:logenergy-differential}
 E_*'+cy\le CR\log(e+R)E_*D(t),
\end{equation}
and consequently
\begin{equation}\label{uni:logenergy}
 \sup_{0\le t\le S}\log E_*(t)
       +\int_0^S\frac{y}{E_*}\,\mathrm dt\le CR\log(e+R).
\end{equation}
\end{lemma}
\begin{proof}
The elliptic estimates, \eqref{uni:coercive}, and the viscosity moment
give
\begin{equation}\label{uni:lowgradient}
\begin{aligned}
 \|\nabla G\|_2+\|\nabla\omega\|_2
       &\le C(\sqrt R A_2+A_3),\\
 |\bar G|&\le CA_3,\qquad
 \|G\|_2\le CR^{b/2}A_1.
\end{aligned}
\end{equation}
Here $\bar G=|\Om|^{-1}\int_\Om\lambda(\rho)\diver u\dd$ and
$\int_\Om\lambda(\rho)\dd\le C$ by \eqref{uni:moment}. We also use
\eqref{uni:comparison} in the form
\begin{equation}\label{uni:energycomparison}
 A_1^2\le C\bigl(A_3^2+\|P-\bar P\|_2^2\bigr),\qquad
 A_3^2\le CE_*.
\end{equation}

Testing momentum by $2\dot u$ gives
\begin{equation}\label{uni:fluxidentity}
\begin{aligned}
 (A_1^2)'+2y={}&
 \int_\Om\left(1+\frac{\rho\lambda'}h\right)\frac{G^2}h\diver u\dd\\
 &-2\int_\Om\left(1+\frac{\rho\lambda'}h\right)G(\diver u)^2\dd
 +2\gamma\int_\Om\frac{GP}h\diver u\dd\\
 &+4\int_\Om G\det\nabla u\dd-\mu\int_\Om\omega^2\diver u\dd
 +2\bar P_t\int_\Om\frac Gh\dd\\
 &-2\int_{\partial\Om}\kappa G u_\tau^2\ds
 +\frac{2\mu}{3}\int_{\partial\Om}(\partial_\tau\alpha)u_\tau^3\ds.
\end{aligned}
\end{equation}
The tangential boundary calculation is the one in \eqref{eng:identity}.
For the normal part, $F=G-\bar P$ and
\[
 2\bar P\int_\Om\diver\dot u\dd
       =-2\bar P\int_{\partial\Om}\kappa u_\tau^2\ds
\]
cancel the mean-pressure boundary terms, leaving the displayed
curvature term with $G$.

Fix $\delta=[4\log(e+R)]^{-1}$. Weighted H\"older and \eqref{eng:GN} give
\[
 \left\|\frac{G^2}h\right\|_2
 \le\frac C{\sqrt\delta}R^{b\delta/2}A_1(\sqrt R A_2+A_3)
             +CA_1^{1-\delta}A_3^{1+\delta}.
\]
Since $\diver u=(G+P-\bar P)/h$, the first three terms of
\eqref{uni:fluxidentity} are bounded by
\[
 CA_3\left(\left\|\frac{G^2}h\right\|_2
          +\left\|\frac{G(P+\bar P)}h\right\|_2\right).
\]
The uniform $L^4$ pressure moment implies
$\|G(P+\bar P)/h\|_2\le C\|G^2/h\|_2^{1/2}$.
Since $\delta^{-1}R^{b\delta}\le C\log(e+R)$, Young's inequality
and \eqref{uni:energycomparison} bound these terms by
\begin{equation}\label{uni:cubic}
 \eps y+C_\eps R\log(e+R)E_*
                    \bigl(A_3^2+\|P-\bar P\|_2^2\bigr).
\end{equation}

For the Jacobian term, choose $p=8+(b+2\gamma+2)\log(e+R)$.
Equations \eqref{eng:GN} and \eqref{uni:lowgradient} yield
\[
 \|G\|_p+\|\omega\|_p
 \le C\sqrt p\,R^{b/p}E_*^{1/p}(\sqrt R A_2+A_3)^{1-2/p}.
\]
Using \eqref{mom:Hodge}, $\diver u=(G+P-\bar P)/h$, and
$\|P-\bar P\|_p\le CR^\gamma$, we obtain
\[
 \|\nabla u\|_p
 \le Cp^{13/2}R^{b/p}E_*^{1/p}(\sqrt R A_2+A_3)^{1-2/p}
                                      +Cp^6R^\gamma.
\]
Interpolation with $\|\nabla u\|_2\le CA_3$ now gives
\begin{equation}\label{uni:Jacobian-interpolation}
\begin{aligned}
 \left|\int_\Om G\det\nabla u\dd\right|
 \le C\sqrt p\,A_3^{2(p-3)/(p-2)}
 \bigl[&R^{b/(p-2)}E_*^{1/(p-2)}(\sqrt R A_2+A_3)\\
       &+R^{b/p+2\gamma/(p-2)}E_*^{1/p}
                              (\sqrt R A_2+A_3)^{1-2/p}\bigr].
\end{aligned}
\end{equation}
The factors $p^{13/(p-2)}$, $p^{12/(p-2)}$, and the displayed extra
powers of $R$ are bounded. When the $A_2$ contributions are
absorbed by Young's inequality, their remaining energy factors satisfy
\[
\begin{aligned}
 E_*^{2/(p-2)}A_3^{4(p-3)/(p-2)}&\le CE_*A_3^2,\\
 E_*^{2/(p+2)}A_3^{4p(p-3)/((p-2)(p+2))}&\le CE_*A_3^2,
\end{aligned}
\]
by \eqref{uni:energycomparison}. The $A_3$ parts satisfy the same
bound. Thus
\begin{equation}\label{uni:Jacobian}
 \left|\int_\Om G\det\nabla u\dd\right|
       \le\eps y+C_\eps R\log(e+R)E_*A_3^2.
\end{equation}

Since $\|\omega\|_2+\|\diver u\|_2\le CA_3$, interpolation gives
\[
 \left|\int_\Om\omega^2\diver u\dd\right|
 \le CA_3^2(\sqrt R A_2+A_3)
 \le\eps y+C_\eps RE_*A_3^2.
\]
The trace inequality and \eqref{uni:lowgradient} similarly give
\[
\begin{aligned}
 \left|\int_{\partial\Om}\kappa G u_\tau^2\ds\right|
 &\le C\|G\|_{H^1}A_3^2
 \le\eps y+C_\eps RE_*A_3^2,\\
 \left|\int_{\partial\Om}(\partial_\tau\alpha)u_\tau^3\ds\right|
 &\le CA_3^3\le CE_*A_3^2.
\end{aligned}
\]
Finally,
\[
 \bar P_t=-\frac{\gamma-1}{|\Om|}\int_\Om(P-\bar P)\diver u\dd,\qquad
 \int_\Om\frac Gh\dd=-\int_\Om\frac{P-\bar P}h\dd,
\]
so
\[
 \left|\bar P_t\int_\Om\frac Gh\dd\right|
 \le CA_3\|P-\bar P\|_2^2
 \le C\bigl(A_3^2+\|P-\bar P\|_4^4\bigr).
\]
Combining these estimates and absorbing the small multiples of $y$
in \eqref{uni:fluxidentity} proves
\eqref{uni:logenergy-differential}. Divide by $E_*$ and use
\eqref{uni:dissipation} and $E_*(0)\le C$ to obtain
\eqref{uni:logenergy}.
\end{proof}

\begin{lemma}[Uniform positive-time acceleration]\label{uni:acceleration-lemma}
There is a fixed $N\ge1$ such that, for $0<\tau<S$,
\begin{equation}\label{uni:acceleration}
 \sup_{\tau\le t\le S}\log(e+y(t))
       +\int_\tau^S\frac{\|\nabla\dot u\|_2^2}{e+y}\,\mathrm dt
       \le C_\tau R^N.
\end{equation}
The constant $C_\tau$ is independent of $S$.
\end{lemma}
\begin{proof}
Testing momentum by $u$ and using the basic kinetic energy gives
\[
 A_3^2\le C A_3\|P-\bar P\|_2
       +A_2\left(\int_\Om\rho|u|^2\dd\right)^{1/2}
       \le CA_3+CA_2.
\]
Together with \eqref{uni:energycomparison}, this implies
\begin{equation}\label{uni:Esquare}
 E_*^2\le C(1+y).
\end{equation}
The fixed-exponent Hodge and flux estimates give
\[
 \|\nabla u\|_4^4
 \le C\left[R^bA_1^2(Ry+A_3^2)+A_3^4+\|P-\bar P\|_4^4\right]
 \le CR^N D(t)(1+y).
\]
In the last step we used $A_1^2+A_3^2\le CD(t)$ and
\eqref{uni:Esquare}.

Use the fixed normal extension from Subsection~\ref{acc:section} and put
\[
 q=\bigl((u\cdot\nabla)n\cdot u\bigr)n,\qquad
 v=\dot u+q,\qquad z=\int_\Om\rho|v|^2\dd.
\]
Then $v\cdot n=0$, and
\begin{equation}\label{uni:qbounds}
\begin{aligned}
 \int_\Om\rho|q|^2\dd&\le CR A_3^4,\\
 \|\nabla q\|_2^2
 &\le C\bigl(A_3^4+A_3^2\|\nabla u\|_4^2\bigr)
       \le CR^N D(t)(1+y).
\end{aligned}
\end{equation}
Define
\[
 D_v=\int_\Om\bigl(h(\diver v)^2+\mu(\curl v)^2\bigr)\dd
            +\mu\int_{\partial\Om}\alpha v_\tau^2\ds,\qquad
 H=1+z+KR E_*^2,
\]
with fixed $K$ sufficiently large. By \eqref{uni:coercive},
\begin{equation}\label{uni:Hcomparison}
 \|v\|_{H^1}^2\le CD_v,\qquad
 1+y\le CH\le CR(1+y).
\end{equation}

Test the commuted equation \eqref{acc:commutation}, corrected by
$\rho\Dt q$, against $v$. The stress remainder in
\eqref{acc:stressremainder} contributes at most
\[
 \eps D_v+CR^{2b}\|\nabla u\|_4^4
                   +C\int_\Om P^2|\nabla u|^2\dd
 \le\eps D_v+CR^N D(t)H,
\]
because
\[
 \int_\Om P^2|\nabla u|^2\dd
 \le C\bigl(A_3^2+\|P-\bar P\|_4^4+\|\nabla u\|_4^4\bigr).
\]
The principal cross terms with $q$ are controlled by
$\eps D_v+CR^b\|\nabla q\|_2^2$.
By \eqref{acc:qbounds}, \eqref{uni:Esquare}, and \eqref{uni:Hcomparison},
\[
\begin{aligned}
 \int_\Om\rho|\Dt q||v|\dd
 &\le C\bigl(\sqrt{Ry}\,A_3+RA_3^3\bigr)D_v^{1/2}\\
 &\le\eps D_v+CR^N D(t)H.
\end{aligned}
\]

Retain the boundary dissipation in \eqref{acc:boundarytraction}.
The trace estimates give
\[
\begin{aligned}
 &\left|\int_{\partial\Om}
       \bigl[\kappa(G-\bar P)u_\tau
                   -\mu(\partial_\tau\alpha)u_\tau^2\bigr]v_\tau\ds\right|\\
 &\qquad\le C\bigl(\sqrt{Ry}\,A_3+A_3+A_3^2\bigr)D_v^{1/2}
       \le\eps D_v+CR^N D(t)H.
\end{aligned}
\]
For the tangential derivative term, the Stokes calculation in
Subsection~\ref{acc:section} gives
\[
 \left|\int_{\partial\Om}\alpha u_\tau(\partial_\tau u_\tau)v_\tau\ds\right|
 \le CD_v^{1/2}\bigl(A_3^2+A_3\|\nabla u\|_4\bigr)
 \le\eps D_v+CR^N D(t)H.
\]
Consequently,
\begin{equation}\label{uni:zdifferential}
 z'+cD_v\le CR^N D(t)H.
\end{equation}
Multiply \eqref{uni:logenergy-differential} by $2KR E_*$ and add
\eqref{uni:zdifferential}. Enlarging the fixed exponent $N$ gives
\begin{equation}\label{uni:Hdifferential}
 H'+cD_v\le CR^N D(t)H.
\end{equation}

By \eqref{uni:logenergy}, there is $s_\tau\in[\tau/2,\tau]$ with
\[
 y(s_\tau)\le\frac2\tau\int_0^\tau y\,\mathrm dt,\qquad
 \log H(s_\tau)\le C_\tau R^N.
\]
Divide \eqref{uni:Hdifferential} by $H$ and use
\eqref{uni:dissipation} to obtain
\[
 \sup_{\tau\le t\le S}\log H(t)
       +\int_\tau^S\frac{D_v}{H}\,\mathrm dt\le C_\tau R^N.
\]
Finally, \eqref{uni:qbounds} and \eqref{uni:Hcomparison} imply
\[
 \frac{\|\nabla\dot u\|_2^2}{e+y}
       \le CR\frac{D_v}{H}+CR^N D(t).
\]
Integration proves \eqref{uni:acceleration}, after enlarging $N$.
\end{proof}

\subsection{The time-uniform density bound}

\begin{lemma}\label{uni:weightedflux-lemma}
For $0<\tau<S$,
\begin{equation}\label{uni:weightedflux}
 \int_\tau^S\frac{\|G\|_\infty^2}{E_*}\,\mathrm dt
                    \le C_\tau R_S^2[\log(e+R_S)]^2.
\end{equation}
Consequently, for $\tau\le s<t\le S$,
\begin{equation}\label{uni:localflux}
 \int_s^t\|G\|_\infty\,\mathrm d\xi
 \le C_\tau R_t\log(e+R_t)
          \left(t-s+\int_s^t A_3^2\,\mathrm d\xi\right)^{1/2}.
\end{equation}
The constants are independent of $S$.
\end{lemma}

\begin{proof}
Fix $R=R_S$. The mass-weighted Poincare inequality gives
$\|\dot u\|_m\le C\sqrt m(A_2+\|\nabla\dot u\|_2)$.
By \eqref{uni:lowgradient}, the mean bound for $G$, and
\eqref{uni:coercive},
\[
 \|G\|_m+\|\omega\|_m\le C\sqrt m(\sqrt R A_2+A_3),\qquad
 \|u\|_m\le C\sqrt m A_3.
\]
Apply \eqref{mom:Hodge} to
$\diver u=(G+P-\bar P)/h$, and then \eqref{mom:fluxelliptic},
to obtain
\begin{equation}\label{uni:highgradient}
 \|\nabla G\|_m
 \le Cm^{15/2}R\bigl(A_2+\|\nabla\dot u\|_2
                                  +A_3+\|P-\bar P\|_m\bigr),
 \qquad m\ge8.
\end{equation}
Retain the pressure norm in this estimate. Its time integral is
controlled by \eqref{uni:pressure-high}.

For this calculation put
\[
 a_0=\frac{A_1}{\sqrt{E_*}},\qquad
 \ell=\frac{A_2+A_3}{\sqrt{E_*}},\qquad
 Z=\frac{\|\nabla\dot u\|_2}{\sqrt{e+y}},\qquad
 \Pi_m=\frac{\|P-\bar P\|_m}{\sqrt{E_*}}.
\]
Equations \eqref{uni:comparison}, \eqref{uni:logenergy},
\eqref{uni:acceleration}, and \eqref{uni:pressure-high} give
\begin{equation}\label{uni:timefactors}
\begin{aligned}
 &0\le a_0\le1,\qquad \int_0^S a_0^2\,\mathrm dt\le C,
       \qquad \int_0^S\ell^2\,\mathrm dt\le CR\log(e+R),\\
 &\int_\tau^S Z^2\,\mathrm dt\le C_\tau R^N,
       \qquad \int_0^S\Pi_m^2\,\mathrm dt\le CmR^N,\qquad
 \sup_{[\tau,S]}\log(e+\ell^2)\le C_\tau R^N.
\end{aligned}
\end{equation}
Enlarge the fixed $N$ to cover these bounds, and choose
\[
 p=8+(N+b+4)\log(e+R),\qquad m=e^p.
\]
Let $\eta,\sigma$ be as in \eqref{cl:interpolation}. Then
$\eta=(1-2/m)\sigma$ and $\sigma-\eta=4/[p(m-2)+2m]$.
The estimates for $G$ and its gradients, inserted in
\eqref{cl:interpolation}, yield
\[
\begin{aligned}
 \frac{\|G-\bar G\|_\infty}{\sqrt{E_*}}
 &\le C\sqrt p\,
 R^{1/2+(b-1)\eta/2+\sigma/2}m^{15\sigma/2}
 a_0^\eta\ell^{1-\eta-\sigma}
       \bigl(\ell+E_*^{-1/2}Z+\ell Z+\Pi_m\bigr)^\sigma\\
 &\le C\sqrt{pR}\,a_0^\eta\ell^{1-\eta-\sigma}
       \bigl(\ell^\sigma+Z^\sigma
                         +\ell^\sigma Z^\sigma+\Pi_m^\sigma\bigr).
\end{aligned}
\]
The extra powers of $R$ and $m$ are bounded. For the first, second,
and fourth summands, apply Holder with time weights
$(\eta,1-\eta)$ or $(\eta,1-\eta-\sigma,\sigma)$ and use
\eqref{uni:timefactors}. For the third, use $a_0\le1$ and
\[
 \sup_{[\tau,S]}\ell^{\sigma-\eta}
 \le\exp\left(\frac{C_\tau R^N}{pm}\right)\le C_\tau,
\]
and then Holder with weights $(1-\sigma,\sigma)$.
Since $R^{N\sigma}$ and $m^\sigma$ are bounded, these four estimates give
\[
 \left(\int_\tau^S
       \frac{\|G-\bar G\|_\infty^2}{E_*}\,\mathrm dt\right)^{1/2}
       \le C_\tau\sqrt{pR}\sqrt{R\log(e+R)}
       \le C_\tau R\log(e+R).
\]
Also $|\bar G|\le CA_3$, so its weighted square integral is bounded.
This proves \eqref{uni:weightedflux}. Apply it on each prefix
$[0,t]$ and use $E_*\le C(1+A_3^2)$ to obtain \eqref{uni:localflux}.
\end{proof}

\begin{proposition}\label{uni:apriori-proposition}
Assume $\alpha\ge0$ and $\alpha\not\equiv0$. With the data of
Proposition~\ref{thm:main}, there is a constant independent of $S$ such that
\begin{equation}\label{uni:apriori}
\begin{aligned}
 &\sup_{0\le t\le S}
       \bigl(\|\rho(t)\|_\infty+\|u(t)\|_{H^1}^2+A_1^2+A_3^2\bigr)
       +\int_0^S A_2^2\,\mathrm dt\\
 &\qquad+\sup_{0<t\le S}\min\{1,t\}A_2^2(t)
       +\int_0^S\min\{1,t\}\|\nabla\dot u\|_2^2\,\mathrm dt\le C.
\end{aligned}
\end{equation}
The dependence is that of Proposition~\ref{thm:main}, with $T$ omitted.
\end{proposition}

\begin{proof}
Use Lemma~\ref{cl:initial} with upper horizon one to fix
$\tau_0,M_0$ independently of $S$. If $S\le\tau_0$,
Proposition~\ref{thm:main} applies. Otherwise take $\tau=\tau_0/2$
and choose $M_*>\max\{2,r_\lambda,M_0\}$, with $\theta(M_*)>0$,
large enough for the absorptions below. Define
\[
 r_j=\theta^{-1}(2^j\theta(M_*)),\qquad t_0=\tau,
\]
and let $t_j$ be the first time the running density maximum reaches
$r_j$. For attained levels put
$\mathcal E_j=\int_{t_{j-1}}^{t_j}A_3^2\,\mathrm dt$; then
$\sum_j\mathcal E_j\le C$ by \eqref{uni:basic}.

On a trajectory ending at density $r_j$, let
$s_j\in[t_{j-1},t_j)$ be its last crossing of $r_{j-1}$.
On $[s_j,t_j]$,
$P-\bar P\ge r_{j-1}^\gamma/2$ after increasing $M_*$.
Since
\[
 \Dt\theta(\rho)+(P-\bar P)=-G,
\]
equation \eqref{uni:localflux} and Young's inequality give
\[
\begin{aligned}
 \frac{\theta(r_j)}2
 &\le C_\tau r_j\log(e+r_j)
                 \sqrt{t_j-s_j+\mathcal E_j}
                -\frac12r_{j-1}^\gamma(t_j-s_j)\\
 &\le C_\tau r_j\log(e+r_j)\sqrt{\mathcal E_j}
        +C_\tau\frac{r_j^2[\log(e+r_j)]^2}{r_{j-1}^\gamma}.
\end{aligned}
\]
Call $j$ good if $\log r_j\le(\gamma+1)\log r_{j-1}/2$.
On a good level the last term is at most
$C_\tau r_j^{2/(\gamma+1)}[\log(e+r_j)]^2
\le\theta(r_j)/4$ for $M_*$ large. Thus
\[
 \mathcal E_j\ge
 c\left(\frac{\theta(r_j)}{r_j\log(e+r_j)}\right)^2
       \ge c[\log(e+r_j)]^{2a-2}\ge\frac{c}{j+1},
\]
where $r_j\le C2^j\theta(M_*)$ and $a>1/2$ were used.
If $j_\ell$ is the $\ell$th bad index, monotonicity gives
\[
 \left(\frac{\gamma+1}{2}\right)^\ell\log M_*
          \le\log r_{j_\ell}\le C(j_\ell+1),\qquad
 \sum_{j\ {\rm bad}}\frac1{j+1}<\infty.
\]
Consequently every attained $L$ satisfies
\[
 C\ge\sum_{j=1}^L\mathcal E_j
       \ge c\sum_{j=1}^L\frac1{j+1}-C.
\]
The harmonic sum diverges, so the number of attained levels, and
hence $R_S$, is bounded independently of $S$.

Equations \eqref{uni:logenergy} and \eqref{uni:comparison} now give
$\sup(A_1^2+A_3^2+\|u\|_{H^1}^2)+\int_0^S y\,\mathrm dt\le C$.
For the acceleration, multiply \eqref{uni:Hdifferential} by
$\sigma(t)=\min\{1,t\}$. Since
$\int_0^{\min\{1,S\}}H\,\mathrm dt\le C$ by
\eqref{uni:Hcomparison}, and $\int_0^S D\,\mathrm dt\le C$,
Gronwall gives
\[
 \sup_{t\le S}\sigma(t)H(t)+\int_0^S\sigma(t)D_v(t)\,\mathrm dt\le C.
\]
The bounds \eqref{uni:qbounds} and \eqref{uni:Hcomparison} transfer
this estimate to $y$ and $\nabla\dot u$, proving \eqref{uni:apriori}.
\end{proof}

\section{Proof of Theorem~\ref{thm:strong}}
\label{sec:strong}

\subsection{Density and velocity derivatives}

For positive approximations, constants may depend on $T$, $q$,
$\|\rho_0\|_{W^{1,q}}$, and the bound of $\lambda'$ on the density
interval controlled by Proposition~\ref{thm:main}.

The low-order estimates and the density-weighted Poincare inequality
for arbitrary vector fields imply
\begin{equation}\label{str:force}
\begin{aligned}
 &\int_0^T\|\rho\dot u\|_q^{(q+1)/q}\,\mathrm dt
   +\int_0^Tt(\|\rho\dot u\|_q^2+\|\dot u\|_{H^1}^2)\,\mathrm dt
       \le C_{q,T},\\
 &\|\nabla F\|_q+\|\nabla\omega\|_q
       \le C_q(\|\rho\dot u\|_q+\|u\|_{W^{1,q}}).
\end{aligned}
\end{equation}
For the first line, interpolate between exponents $2$ and $q^2$,
using
\[
 \|\rho\dot u\|_q\le C\left(A_2+
  A_2^{2(q-1)/(q^2-2)}
       \|\nabla\dot u\|_2^{q(q-2)/(q^2-2)}\right).
\]
After the power $(q+1)/q$ and time Holder, the remaining singular
time exponent is
$q(q-2)(q+1)/((q^2-2)(q-1))<1$.
The weighted assertion follows directly from the weighted acceleration.

At exponent two, \eqref{eng:ellipticlow} gives
$\|F\|_{H^1}+\|\omega\|_{H^1}\le C_T(1+A_2)$.
The div--curl estimate with its lower order term then gives, for each
fixed finite $p$,
\begin{equation}\label{str:uWp}
 \|u\|_{W^{1,p}}\le C_{p,T}(1+A_2).
\end{equation}
In particular \eqref{str:force} gives
$\|\nabla F\|_q+\|\nabla\omega\|_q
 \le C_{q,T}(1+A_2+\|\rho\dot u\|_q)$.

We use the bounded-domain logarithmic div--curl estimate
\begin{equation}\label{str:BKM}
 \|\nabla u\|_\infty
 \le C_q\bigl(1+\|\diver u\|_\infty+\|\omega\|_\infty
                 +\|u\|_{H^1}\bigr)
                      \log(e+\|u\|_{W^{2,q}}).
\end{equation}
Write
$u=\nabla\phi+\nabla^\perp\chi+h_0$ using the scalar homogeneous
Neumann and Dirichlet potentials from Lemma~\ref{mom:potential}.
The finite-dimensional harmonic part satisfies
$\|\nabla h_0\|_\infty+\|h_0\|_{W^{2,q}}\le C\|u\|_2$.
The scalar endpoint Hessian estimate gives
\[
 [\nabla^2\phi+\nabla\nabla^\perp\chi]_{\mathrm{BMO}}
       \le C(\|\diver u\|_\infty+\|\omega\|_\infty).
\]
For the scalar endpoint estimate, see
\cite[Appendix~A, Lemma~A.3]{HilhorstRoux}. A local proof follows.
For $\Delta w=f$ with either homogeneous boundary condition, use the
fixed $W^{2,4}$ estimate. After fixed localization and odd Dirichlet
or even Neumann reflection in normal coordinates, the equation is
\[
 A(y):D^2W+B(y)\cdot\nabla W=G,
 \qquad
 \|D^2W\|_4+\|\nabla W\|_\infty+\|G\|_\infty\le C\|f\|_\infty,
\]
where $A$ is Lipschitz and uniformly elliptic, and $B$ is bounded.
On a small reflected ball $B_r$, freeze $A_0=A(y_0)$ and solve
\[
 A_0:D^2V=G-B\cdot\nabla W+(A_0-A):D^2W,
       \qquad V|_{\partial B_r}=0.
\]
The root mean square of the right hand side is at most
\[
 C\|f\|_\infty+Cr\left(\frac1{|B_r|}\int_{B_r}|D^2W|^2\,\mathrm dy\right)^{1/2}
       \le C\|f\|_\infty+Cr^{1/2}\|D^2W\|_4\le C\|f\|_\infty.
\]
The scale-invariant $L^2$ elliptic estimate therefore gives
$|B_r|^{-1}\int_{B_r}|D^2V|^2\,\mathrm dy\le C\|f\|_\infty^2$. Each component of
$D^2(W-V)$ is $A_0$-harmonic. For
$\Psi(r)=|B_r|^{-1}\int_{B_r}|D^2W-(D^2W)_{B_r}|^2\,\mathrm dy$,
where $(D^2W)_{B_r}$ denotes the ball average, homogeneous
oscillation decay and this decomposition imply
\[
 \Psi(\vartheta r)\le C\vartheta^2\Psi(r)
                                  +C\vartheta^{-2}\|f\|_\infty^2.
\]
Choose a fixed $\vartheta$ with $C\vartheta^2<1/2$.
Iteration and the fixed-radius $L^4$ bound give
$\sup_r\Psi(r)\le C\|f\|_\infty^2$.
Returning to physical coordinates adds bounded gradient terms and
coefficient oscillations at most $Cr^{1/2}\|D^2W\|_4\le C\|f\|_\infty$.
A finite chart cover and the fixed $L^4$ estimate on larger balls
give the domain BMO bound.
Telescope averages of $\nabla u-\nabla h_0$ on concentric balls
intersected with the domain. Each scale costs at most
$C(\|\diver u\|_\infty+\|\omega\|_\infty)$, the largest average
is bounded by $C\|u\|_{H^1}$, and the Morrey error is at most
$C2^{-L(1-2/q)}\|u\|_{W^{2,q}}$.
Choose $L\asymp\log(e+\|u\|_{W^{2,q}})$ to prove
\eqref{str:BKM}.

Put $\Phi_i=(2\mu+\lambda(\rho))\partial_i\rho$. Direct substitution
of the divergence relation in the differentiated continuity equation
gives the cancellation
\begin{equation}\label{str:Phi}
 \Dt\Phi_i+(2\mu+\lambda)\partial_i u_j\partial_j\rho
        +\Phi_i\diver u+\rho\partial_iF+\rho\partial_iP=0.
\end{equation}
The $\rho\lambda'\partial_i\rho\diver u$ terms cancel exactly.
Testing \eqref{str:Phi} gives
\[
 (\|\Phi\|_q)'
 \le C(1+\|\nabla u\|_\infty)\|\Phi\|_q
                  +C_q\|\nabla F\|_q.
\]
The fixed lower $L^2$ flux norms and Gagliardo--Nirenberg imply
\[
 \|F\|_\infty+\|\omega\|_\infty+\|\diver u\|_\infty
 \le C_{q,T}\left(1+A_2^{q/(2(q-1))}
                  +\|\rho\dot u\|_q^{q/(2(q-1))}\right).
\]
Differentiate $\diver u=(F+P)/(2\mu+\lambda)$ once and use the local
bound for $\lambda'$ to obtain
\[
\begin{aligned}
 \|u\|_{W^{2,q}}
 \le C_{q,T}\bigl[&\bigl(1+A_2^{q/[2(q-1)]}
              +\|\rho\dot u\|_q^{q/[2(q-1)]}\bigr)\|\Phi\|_q\\
                 &+1+A_2+\|\rho\dot u\|_q\bigr].
\end{aligned}
\]
Since $q/[2(q-1)]<1$, using
$(1+s^{q/[2(q-1)]})\log(e+s)\le C_q(1+s)$ in \eqref{str:BKM} gives
\[
 \|\nabla u\|_\infty
 \le C_{q,T}(1+A_2+\|\rho\dot u\|_q)\log(e+\|\Phi\|_q).
\]
The coefficient is integrable by \eqref{str:force}. Logarithmic
Gronwall proves
\begin{equation}\label{str:densitygradient}
 \sup_{t\le T}\|\rho(t)\|_{W^{1,q}}
            +\int_0^T\|\nabla u\|_\infty\,\mathrm dt\le C_{q,T}.
\end{equation}
The preceding elliptic estimate now yields both the unweighted
$L^{(q+1)/q}_tW^{2,q}$ bound and the weighted $L^2_tW^{2,q}$ bound.
At exponent two interpolation absorbs
$\|\nabla\rho\|_q\|\diver u\|_{2q/(q-2)}$, giving
\[
 \|u\|_{H^2}\le C_{q,T}(1+A_2).
\]
The estimates for $u_t=\dot u-u\cdot\nabla u$ follow from
\[
 \|u\cdot\nabla u\|_2^2\le C_{q,T}(1+\|u\|_{H^2}^2),
 \quad
 \|\nabla(u\cdot\nabla u)\|_2^2
 \le C_{q,T}(1+A_2^2+\|u\|_{W^{2,q}}^2).
\]
Finally, the continuity equation gives $\sup\|\rho_t\|_2\le C_{q,T}$.

\subsection{Positive local solutions and continuation}

For the positive-density construction, let the bulk coefficient be
smooth. Local constants may depend on the density lower bound and
coefficient derivatives on the attained density interval.

\begin{lemma}[Linearized equation with transported coefficients]
\label{loc:linear}
Let $\rho\in C([0,t_*];H^2)$, $0<c\le\rho\le C$,
$\sup_{t\le t_*}\|\rho(t)\|_{H^2}\le L$,
and $\rho_t+\diver(\rho v)=0$. Assume that $v$ is tangent and
\[
 v\in C H^2\cap L^2H^3,\qquad
 v_t\in L^\infty L^2\cap L^2H^1,\qquad
 \rho_t\in L^\infty H^1,\qquad \|v\|_{L^\infty_{t,x}}\le V.
\]
Define
\[
 \mathcal L_\rho u=-\mu\Delta u
             -\nabla((\mu+\lambda(\rho))\diver u).
\]
For $u_0\in H^2$ satisfying both Navier trace conditions, the problem
\begin{equation}\label{loc:equation}
 \rho(u_t+v\cdot\nabla u)+\mathcal L_\rho u=-\nabla P
\end{equation}
with \eqref{eq:slip} has a unique solution
\[
 u\in C H^2\cap L^2H^3,\qquad
 u_t\in L^\infty L^2\cap L^2H^1.
\]
Set $Z=\int_\Om\rho(|u|^2+|u_t|^2)\dd$. Then
\begin{align}
 \|u\|_{H^2}&\le C_{L,V}(1+\|u_t\|_2+\|u\|_2),
                                                   \label{loc:H2}\\
 Z'+c_1(\|\nabla u\|_2^2+\|\nabla u_t\|_2^2)
 &\le C_{L,V}\bigl(1+\|\rho_t\|_{H^1}^2+\|v_t\|_2^{4/3}\bigr)(1+Z),
                                                   \label{loc:energy}\\
 \|u\|_{H^3}&\le C_L\left[\|u_t\|_{H^1}
              +(1+\|v\|_{H^2})\|u\|_{H^2}+1\right].
                                                   \label{loc:H3}
\end{align}
The constants also depend on $c,C$, the geometry and boundary coefficient,
and the coefficient on $[c,C]$.
\end{lemma}

\begin{proof}
For tangent fields the elliptic operator has the variational form
\[
 a_\rho(w,z)=\int_\Om\bigl((2\mu+\lambda(\rho))\diver w\diver z
                         +\mu\curl w\curl z\bigr)\dd
                    +\mu\int_{\partial\Om}\alpha w_\tau z_\tau\ds.
\]
It satisfies
\[
 \|w\|_{H^1}^2\le C\bigl(a_\rho(w,w)+\|w\|_2^2\bigr).
\]

The flux and vorticity equations, $(2\mu+\lambda(\rho))^{-1}\in H^2$,
and the div--curl estimate give
\[
 \|u\|_{H^2}\le C_L(1+\|\rho(u_t+v\cdot\nabla u)\|_2+\|u\|_{H^1}).
\]
Using $\|v\|_\infty\le V$ and interpolating between $L^2$ and $H^2$
proves \eqref{loc:H2}.

Differentiate \eqref{loc:equation} in time. In variational form,
\begin{equation}\label{loc:Uequation}
\begin{aligned}
 \rho(u_{tt}+v\cdot\nabla u_t)+\mathcal L_\rho u_t
 ={}&-\rho_tu_t-\rho_tv\cdot\nabla u-\rho v_t\cdot\nabla u\\
    &+\nabla(\lambda'(\rho)\rho_t\diver u)-\nabla P_t.
\end{aligned}
\end{equation}
The fixed boundary operator gives the Navier condition for $u_t$
in variational form. Testing by $u_t$ gives the pressure and viscosity
terms
\[
 \int_\Om P_t\diver u_t\dd,\qquad
       -\int_\Om\lambda'(\rho)\rho_t\diver u\diver u_t\dd.
\]
They are controlled by
\[
 \|P_t\|_2\le C_L\|\rho_t\|_2,\qquad
 \|\lambda'(\rho)\rho_t\diver u\|_2
            \le C_L\|\rho_t\|_{H^1}\|u\|_{H^2}.
\]
The transport equation cancels the mass derivative and advection.
The remaining terms satisfy
\[
\begin{aligned}
 &\left|\int_\Om\rho_t|u_t|^2\dd\right|
  +\left|\int_\Om\rho_tv\cdot\nabla u\,u_t\dd\right|\\
 &\qquad\le\eps\|\nabla u_t\|_2^2
  +C_{\eps,L,V}(1+\|\rho_t\|_{H^1}^2)
                         (\|u\|_{H^2}^2+\|u_t\|_2^2),\\
 \left|\int_\Om\rho v_t\cdot\nabla u\,u_t\dd\right|
 &\le C\|v_t\|_2\|u\|_{H^2}\|u_t\|_2^{1/2}\|u_t\|_{H^1}^{1/2}\\
 &\le\eps\|\nabla u_t\|_2^2
       +C_{\eps,L,V}(1+\|v_t\|_2^{4/3})
                         (\|u\|_{H^2}^2+\|u_t\|_2^2).
\end{aligned}
\]
Testing \eqref{loc:equation} also by $u$, adding, and using
\eqref{loc:H2} proves \eqref{loc:energy}.

For \eqref{loc:H3}, we have
\[
 \|\rho(u_t+v\cdot\nabla u)\|_{H^1}
       \le C_L\bigl(\|u_t\|_{H^1}+\|v\|_{H^2}\|u\|_{H^2}\bigr).
\]
The Dirichlet equation
\[
 \mu\Delta\omega=\curl\bigl(\rho(u_t+v\cdot\nabla u)\bigr)
\]
with boundary value $-\alpha u_\tau\in H^{3/2}(\partial\Om)$
bounds $\omega$ in $H^2$. The identity
\[
 \nabla F=\rho(u_t+v\cdot\nabla u)-\mu\nabla^\perp\omega
\]
and the mean of
$F=(2\mu+\lambda)\diver u-P$ then bound $F$ in $H^2$.
Use the $H^2$ algebra property in two dimensions to estimate
$(F+P)/(2\mu+\lambda(\rho))$ in $H^2$, and then the higher
div--curl estimate, including its harmonic part, to obtain \eqref{loc:H3}.

Construct $u$ at times $t_n=nh$ by the implicit variational problem
\[
 \rho^n\frac{u^n-u^{n-1}}h+\mathcal L_{\rho^n}u^n
          +\rho^nv^n\cdot\nabla u^n=-\nabla P^n.
\]
For small $h$ the positive mass term absorbs the negative part of
the advection form. Lax--Milgram and elliptic regularity give
$u^n\in H^2$ with the Navier conditions. Set
$U^n=(u^n-u^{n-1})/h$, and define the initial $L^2$ value by
\begin{equation}\label{loc:U0}
 U^0=\rho_0^{-1}\left[-\mathcal L_{\rho_0}u_0
          -\rho_0v_0\cdot\nabla u_0-\nabla P_0\right].
\end{equation}
Subtract the preceding step, using its coefficients as the principal
coefficients, to obtain
\[
\begin{aligned}
 \rho^{n-1}\frac{U^n-U^{n-1}}h
 +\mathcal L_{\rho^{n-1}}U^n
 +\rho^{n-1}v^{n-1}\cdot\nabla U^n
 ={}&-\delta_h\rho^n U^n
       -\delta_h(\rho v)^n\cdot\nabla u^n\\
 &+\nabla(\delta_h\lambda^n\diver u^n)-\nabla\delta_hP^n.
\end{aligned}
\]
For $n=1$, subtract \eqref{loc:U0}; the principal operators still act
on $U^1$, so only $U^0\in L^2$ is needed. Set
$Y_n=\int_\Om\rho^n|U^n|^2\dd$ and
$Z_n=\int_\Om\rho^n(|u^n|^2+|U^n|^2)\dd$. The mass identities are
\[
\begin{aligned}
 \int_\Om\rho^{n-1}\delta_hU^n\cdot U^n\dd
 &=\frac{Y_n-Y_{n-1}}{2h}
  +\frac1{2h}\int_\Om\rho^{n-1}|U^n-U^{n-1}|^2\dd\\
 &\quad-\frac12\int_\Om\delta_h\rho^n|U^n|^2\dd,\\
 \int_\Om\rho^n U^n\cdot u^n\dd
 &=\frac{\int_\Om\rho^n|u^n|^2\dd-\int_\Om\rho^{n-1}|u^{n-1}|^2\dd}{2h}\\
 &\quad+\frac h2\int_\Om\rho^n|U^n|^2\dd
   -\frac12\int_\Om\delta_h\rho^n|u^{n-1}|^2\dd.
\end{aligned}
\]
The previous-step elliptic bound gives
\[
 \|u^{n-1}\|_{H^1}^2
       \le C_{L,V}(1+\|U^{n-1}\|_2^2+\|u^{n-1}\|_2^2)
       \le C_{L,V}(1+Z_{n-1}).
\]
For $n=1$ this follows from \eqref{loc:U0}. The advection terms
are $\frac12\int_\Om\rho_t(t_n)|u^n|^2\dd$ and
$\frac12\int_\Om\rho_t(t_{n-1})|U^n|^2\dd$, respectively, and are
estimated by $L^4$ interpolation. Here
$\rho_t=-\diver(\rho v)\in C H^1$. Moreover,
$\delta_h\rho^n,\delta_hv^n$ are time averages and
\[
 \delta_h(\rho v)^n=(\delta_h\rho^n)v^n+\rho^{n-1}\delta_hv^n.
\]
The mean value formula bounds $\delta_h\lambda^n,\delta_hP^n$
by $\delta_h\rho^n$. Thus, with
\[
 C_n=C_{L,V}\bigl(1+\|\rho_t\|_{L^\infty H^1}^2
                              +\|\delta_hv^n\|_2^{4/3}\bigr),
\]
the tested equations yield
\[
 \frac{Z_n-Z_{n-1}}h
    +c_1(\|\nabla u^n\|_2^2+\|\nabla U^n\|_2^2)
       \le C_n(1+Z_n+Z_{n-1}).
\]
Jensen's inequality gives
$h\sum_n\|\delta_hv^n\|_2^{4/3}
\le\int_0^{t_*}\|v_t\|_2^{4/3}\,\mathrm dt$, and the density node terms are
bounded by $2t_*\|\rho_t\|_{L^\infty H^1}^2$.
Since $v_t\in L^\infty L^2$, the $C_n$ are uniformly bounded.
For $h\sup_n C_n$ small, discrete Gronwall and \eqref{loc:H2}
give uniform $L^\infty H^2$, $L^\infty L^2$ bounds for $u^n,U^n$
and an $L^2H^1$ bound for $U^n$; \eqref{loc:H3} then gives the
$L^2H^3$ bound.

Use the affine interpolation $u_h$ for the time derivative and the
backward constant interpolation $\bar u_h$ for the $H^3$ bound.
Since
\[
 \|u_h-\bar u_h\|_{L^2H^1}\le h\|U_h\|_{L^2H^1}\longrightarrow0,
\]
weak compactness and coefficient continuity give a common limit
satisfying the variational equation, with $u_t\in L^2H^1$ and
$u\in L^2H^3$.
The Hilbert triple then gives $u\in C H^2$. Time differences of the
limiting equation, or the same discrete energy limit, justify
\eqref{loc:Uequation}--\eqref{loc:energy}. Linear uniqueness follows
from the weighted $L^2$ energy identity.
\end{proof}

\begin{lemma}[Positive local solution]\label{loc:existence}
For positive $\rho_0\in H^2$ and $u_0\in H^2$ satisfying the Navier
trace conditions, a smooth nonnegative bulk coefficient gives a unique
local solution
\[
 \rho,u\in C H^2,\qquad u\in L^2H^3,\qquad
 u_t\in L^\infty L^2\cap L^2H^1.
\]
The local time may depend on the positive density lower bound, the
$H^2$ data, and the coefficient on the relevant positive interval.
\end{lemma}

\begin{proof}
For tangent $v$ with $v(0)=u_0$, transport $\rho_0$ along $v$ and
apply Lemma~\ref{loc:linear}. This defines the map $v\mapsto u$.

Fix the density bounds and $L>\|\rho_0\|_{H^2}$ using the initial
data, and put $V=\|u_0\|_\infty+1$. For a radius $K$, use the ball
with $v(0)=u_0$, the Navier conditions, and
\begin{equation}\label{loc:ball}
 \|v\|_{L^\infty H^2}+\|v_t\|_{L^\infty L^2}
       +\|v_t\|_{L^2H^1}+\|v\|_{L^2H^3}\le K.
\end{equation}
The transport estimate and Sobolev embedding give
\[
 \|\rho(t)\|_{H^2}\le\|\rho_0\|_{H^2}
             \exp(C\sqrt{t_*}K),\qquad
 \|\rho_t\|_{H^1}\le C_LK.
\]
The flow formula gives the analogous positive lower and upper bounds.
Moreover,
\[
 \|v-u_0\|_\infty
 \le C\|v-u_0\|_{H^1}^{1/2}\|v-u_0\|_{H^2}^{1/2}
 \le Ct_*^{1/4}K^{1/2}(K+\|u_0\|_{H^2})^{1/2}.
\]
After fixing $K$, choose $t_*$ so that the density bounds and
$\|v\|_\infty\le V$ hold for the fixed $L,V$.

The initial value \eqref{loc:U0} is fixed by the data since $v_0=u_0$.
Equation~\eqref{loc:energy} has integrated coefficient at most
$C_{L,V}t_*(1+K^2+K^{4/3})$. Hence, when this is small,
\[
 \|u\|_{L^\infty H^2}+\|u_t\|_{L^\infty L^2}
                 +\|u_t\|_{L^2H^1}\le C_{\rm data},
 \qquad
 \|u\|_{L^2H^3}\le C_{\rm data}(1+\sqrt{t_*}K).
\]
The Navier Hodge semigroup applied to $u_0$ has the norms in
\eqref{loc:ball} bounded by $C\|u_0\|_{H^2}$ for $t_*\le1$.
Choose $K$ larger than four times the sum of this bound and the
displayed data constants, then choose $t_*$ to satisfy the smallness
requirements. The ball is nonempty and the map preserves it.

For two input velocities write $\delta v=v_1-v_2$,
$r=\rho_1-\rho_2$, and $w=u_1-u_2$. The transported density
difference satisfies
\begin{equation}\label{loc:contractrho}
 \sup_{s\le t_*}\|r(s)\|_2^2
       \le C_K t_*\int_0^{t_*}\|\delta v\|_{H^1}^2\,\mathrm dt.
\end{equation}
Indeed, $\nabla\rho_i\in L^\infty L^4$ and
$\int_0^{t_*}\|v_i\|_{W^{1,\infty}}\,\mathrm dt\le C\sqrt{t_*}K$.
The output difference equation is
\[
\begin{aligned}
 \rho_1(w_t+v_1\cdot\nabla w)+\mathcal L_{\rho_1}w
 ={}&-r(u_{2,t}+v_2\cdot\nabla u_2)
                  -\rho_1\delta v\cdot\nabla u_2\\
   &+\nabla((\lambda(\rho_1)-\lambda(\rho_2))\diver u_2)
                         -\nabla(P_1-P_2).
\end{aligned}
\]
Testing by $w$, use
$u_{2,t}+v_2\cdot\nabla u_2\in L^2H^1$, the local Lipschitz bound
for $\lambda$, and
\[
 \left|\int_\Om\rho_1\delta v\cdot\nabla u_2\,w\dd\right|
 \le C\|\nabla u_2\|_\infty(\|\delta v\|_2^2+\|w\|_2^2).
\]
Positivity of $\rho_1$ and the div--curl estimate give, with
$D(w)=\sup_{t\le t_*}\|w(t)\|_2^2+\int_0^{t_*}\|w\|_{H^1}^2\,\mathrm dt$,
\[
 D(w)\le C_K(\sqrt{t_*}+t_*)D(\delta v).
\]
Here \eqref{loc:contractrho} controls the density sources, and
$\int_0^{t_*}\|\nabla u_2\|_\infty\,\mathrm dt\le C\sqrt{t_*}K$ controls the input
velocity source. For $t_*$ small the map contracts in the distance
$D(v_1-v_2)^{1/2}$. Weak compactness in \eqref{loc:ball} makes the
ball closed in this metric and preserves its initial value and boundary
conditions. Banach iteration gives a solution, and the difference
estimate gives uniqueness. The transport equation and the Hilbert
triple give $\rho,u\in C H^2$.
\end{proof}

To apply Proposition~\ref{thm:main}, note that
\[
 \rho_t\in L^\infty H^1,\quad
 F,\omega\in L^\infty H^1\cap H^1(0,t_*;L^2),\quad
 \dot u\in L^2H^1.
\]
Smooth positive-part truncations and the transport formula justify
the density moments. Time differences and $H^1$ traces give the flux
energy identity. For the acceleration identity, take time differences
of the variational momentum equation and test by the Steklov-regularized
corrected acceleration from Subsection~\ref{acc:section}. The boundary
terms are interpreted variationally and yield \eqref{acc:boundarytraction};
its factors have the required traces since $u\in L^\infty H^2$ and
the corrected acceleration belongs to $L^2H^1$.
The initial $L^2$ acceleration is finite for each positive
approximation, so the time-weighted initial terms vanish.

For positive smooth data, suppose the maximal $H^2$ time $T^*$ were
finite. The a priori estimates and \eqref{str:densitygradient} hold
uniformly on $[0,S]$, $S<T^*$. The flow formula gives
\[
 \inf_{t<T^*,x}\rho(t,x)
 \ge\inf\rho_0\exp\left(-\int_0^{T^*}\|\diver u\|_\infty\,\mathrm dt\right)>0.
\]
In the remainder of this continuation step constants may depend on
that lower bound and the higher initial norms. The corrected
acceleration inequality, without a time weight, gives
$\sup A_2^2+\int_0^S\|\dot u\|_{H^1}^2\,\mathrm dt\le C_*$.
Thus $\sup\|u\|_{H^2}\le C_*$ and
\[
 \int_0^S(\|\nabla^2F\|_2^2+\|\nabla^2\omega\|_2^2)\,\mathrm dt\le C_*.
\]
The vorticity boundary lift here is controlled by $u$ in $H^2$.

To control the second density derivatives, we claim that
\[
 (\|\nabla\Phi\|_2)'
 \le C_*(1+\|\nabla u\|_\infty)(1+\|\nabla\Phi\|_2)
        +C_*(\|\nabla^2F\|_2+\|\nabla^2\omega\|_2).
\]
Indeed, $\nabla\Phi=(2\mu+\lambda)\nabla^2\rho
+\lambda'\nabla\rho\otimes\nabla\rho$ and $q>2$ imply
\[
 \|\nabla^2\rho\|_2+\|\nabla\rho\|_4^2
                      \le C_*(1+\|\nabla\Phi\|_2).
\]
For $2<q<4$, the interpolation power of $\|\nabla^2\rho\|_2$
in $\|\nabla\rho\|_4^2$ is $2-q/2<1$. Also,
\[
 \|\nabla\rho\nabla F\|_2
                 \le C_*(1+\|\nabla^2F\|_2^{2/q}),\qquad
 \|\nabla\rho\nabla^2u\|_2
                     \le C_*(1+\|u\|_{H^3}^{2/q}).
\]
The second derivative of $\diver u=(F+P)/(2\mu+\lambda)$ and
the div--curl estimate give
\[
 \|u\|_{H^3}\le C_*\bigl[1+\|\nabla^2F\|_2+\|\nabla^2\omega\|_2
 +(1+\|\nabla u\|_\infty)(1+\|\nabla\Phi\|_2)
 +\|u\|_{H^3}^{2/q}\bigr].
\]
Absorbing the last term and differentiating \eqref{str:Phi} proves
the claim. Gronwall gives $\sup\|\rho\|_{H^2}\le C_*$.
Now use $\|\nabla u\|_\infty
\le\delta\|u\|_{H^3}+C_\delta\|u\|_{H^2}$ in the preceding
elliptic inequality and absorb to obtain
$\|u\|_{H^3}\le C_*(1+\|\nabla^2F\|_2+\|\nabla^2\omega\|_2)$.
Hence $u\in L^2H^3$.
In the twice-differentiated continuity equation, the source terms
$\rho\nabla^2\diver u$, $\nabla\rho\nabla^2u$ and
$\nabla u\nabla^2\rho$ belong to $L^1L^2$.
Transport by the integrably Lipschitz flow therefore gives the
terminal $H^2$ strong density trace. Also
$u_t=\dot u-u\cdot\nabla u\in L^2H^1$; the Hilbert triple gives
the terminal $H^2$ velocity trace. They retain positive density and the
Navier boundary condition. Lemma~\ref{loc:existence} therefore restarts
the local solution from these terminal data.
This contradicts maximality.

\subsection{Passage to vacuum}

For $0<\varepsilon<1$, choose a nonnegative smooth kernel $\chi$ supported in $(-1,1)$ with
integral one, and define
\[
 \lambda_\varepsilon(r)=\int_{-1}^1\chi(z)
                   \lambda(r+\varepsilon(2+z))\,\mathrm dz.
\]
These are smooth and nonnegative, converge locally in $C^1$ on
$[0,\infty)$, and satisfy
\[
 |r\lambda_\varepsilon'(r)|\le C_\lambda(1+\lambda_\varepsilon(r)).
\]
For large $r$ the shifted arguments lie in $[r,2r]$. The logarithmic
derivative condition implies
$\lambda_\varepsilon(r)\ge c\lambda(r)-C$, and hence
\[
 \theta_\varepsilon(r)\ge c\theta(r)-C\log r-C.
\]
The growth condition is preserved with common constants after
increasing a fixed threshold. Local bounds of
$\lambda_\varepsilon'$ are uniform.

Take positive smooth approximations of $\rho_0$ converging in
$W^{1,q}$, with common mass and common $L^\infty$ bound. For the velocity
use $u_{0,\varepsilon}=e^{-\varepsilon L_\alpha}u_0$, where $L_\alpha$
is the nonnegative Navier Hodge operator associated with the form
\[
 \int_\Om(\diver v\diver w+\curl v\curl w)\dd
                   +\int_{\partial\Om}\alpha v_\tau w_\tau\ds
\]
on tangent $H^1$ fields. The form norm plus the $L^2$ norm is equivalent
to $H^1$, so these approximations converge in $H^1$, satisfy the
Navier conditions, and preserve any harmonic component when $\alpha=0$.
The positive construction gives global approximate solutions;
its local constants may depend on $\varepsilon$.
The common structural constants and local $C^1$ bounds for
$\lambda_\varepsilon$ make \eqref{eq:main}, \eqref{str:force}, and
\eqref{str:densitygradient} uniform on each finite interval.
For $\alpha\not\equiv0$, Proposition~\ref{uni:apriori-proposition}
also gives \eqref{uni:apriori} with a constant independent of that
interval and of the approximation parameter.

Compactness gives, for $0<\eta<1-2/q$,
\[
 \rho_\varepsilon\to\rho\text{ in }C([0,T];C^\eta(\overline\Om)),
 \qquad u_\varepsilon\to u\text{ in }L^2(0,T;H^1).
\]
For the velocity first work on $[\tau,T]$ using the weighted time
derivative and $W^{2,q}$ bounds; the common $L^\infty H^1$ norm
controls the omitted interval as $\tau\downarrow0$.
The momentum is bounded in $L^\infty H^1$ and its time derivative
in $L^\infty H^{-1}$: test the conservative equation against
$H^1_0$ fields, so only the interior divergence of the uniformly
$L^2$ stress and convective flux is needed. Thus
$\rho_\varepsilon u_\varepsilon\to\rho u$ in $C([0,T];L^2)$.
Uniform convergence of $\lambda_\varepsilon(\rho_\varepsilon)$ identifies
the viscous stress. The strong velocity convergence identifies the
convective term. On each $[\tau,T]$, the weak $L^2W^{2,q}$ convergence
and the continuous first-derivative trace map pass
$\curl u_\varepsilon=-\alpha u_{\varepsilon,\tau}$ to the limit.
The normal trace follows already from the $H^1$ convergence.
Weak compactness of derivatives
gives the equations almost everywhere at positive times.

The differentiated transport formula, using
$\nabla u\in L^1L^\infty$ and $\rho\nabla\diver u\in L^1L^q$,
gives $\rho\in C W^{1,q}$. For $\delta>0$,
$\rho u/\sqrt{\rho+\delta}\in C L^2$, and its distance from
$\sqrt\rho\,u$ is at most $\sqrt\delta\|u\|_2$, uniformly in time.
This proves the stated weighted velocity trace. The positive-time
time derivative bounds identify $u_t$ and $\dot u$; lower
semicontinuity followed by $\tau\downarrow0$ retains all weights.
A diagonal sequence over integer time intervals gives one global
solution.

\subsection{Uniqueness}

Let $r=\rho_1-\rho_2$ and $w=u_1-u_2$ for two solutions with the
same initial density and momentum. Define
\[
\begin{aligned}
 Y&=\int_\Om\rho_1|w|^2\dd,\\
 D_w&=\int_\Om((2\mu+\lambda(\rho_1))(\diver w)^2+\mu(\curl w)^2)\dd
                      +\mu\int_{\partial\Om}\alpha w_\tau^2\ds.
\end{aligned}
\]
The weighted div--curl inequality gives
$\|w\|_{H^1}^2\le C(D_w+Y)$.
The density difference equation yields
\begin{equation}\label{str:diffdensity}
 \|r(t)\|_2^2\le Ct\int_0^t(D_w+Y)\,\mathrm ds.
\end{equation}
Here $\nabla\rho_1\in L^\infty L^q$ and
$\nabla u_2\in L^1L^\infty$.

On the common bounded density range, $P$ and $\lambda$ are Lipschitz.
Testing the momentum difference by $w$ and using the weighted
div--curl inequality and the Navier condition gives
\begin{equation}\label{str:diffenergy}
 Y'+cD_w
 \le C(1+\|\nabla u_2\|_\infty)Y
   +C(1+\|\nabla u_2\|_\infty^2+\|\dot u_2\|_q^2)\|r\|_2^2.
\end{equation}
The viscosity difference contributes at most
$\eps D_w+C_\eps\|\nabla u_2\|_\infty^2\|r\|_2^2$.

The strong and weighted estimates imply
\[
 \int_0^T\left[1+\|\nabla u_2\|_\infty+
  t(1+\|\nabla u_2\|_\infty^2+\|\dot u_2\|_q^2)\right]\,\mathrm dt<\infty.
\]
The common initial momentum and weighted velocity traces give
$Y(0+)=0$. Set
$\mathcal Y(t)=\sup_{s\le t}Y(s)+\int_0^tD_w\,\mathrm ds$.
Substitute \eqref{str:diffdensity} into \eqref{str:diffenergy}, integrate
from a positive time, and let it tend to zero. The displayed
integrable coefficient and Gronwall give $\mathcal Y=0$.
Thus $\rho_1=\rho_2$ and $u_1=u_2$. This completes
the existence and uniqueness assertion. Lower semicontinuity in
\eqref{uni:apriori}, followed by $S\to\infty$, gives
\eqref{str:uniform} and proves Theorem~\ref{thm:strong}.
For general $\alpha\ge0$, the same construction and uniqueness
argument use Proposition~\ref{thm:main} on each finite interval;
this proves Remark~\ref{rem:general-slip}.

\section{Large-time behavior}\label{sec:longtime}

We prove Theorem~\ref{thm:longtime}. The estimates below are first
obtained for the positive approximations and pass to the strong
solution. Constants are independent of time and of the approximation
parameter unless a positive density lower bound is explicitly assumed.

\subsection{Exponential decay}

The uniform density bound in \eqref{uni:apriori} fixes $R$ in
Lemma~\ref{uni:conditional-decay-lemma}. Its proof gives exponential
decay of the relative energy $\mathcal H$ in \eqref{uni:relative-energy}.
Since $\mathcal H'+A_3^2=0$ and the density is bounded,
\[
 \mathcal H\asymp\|\sqrt\rho\,u\|_2^2+\|\rho-\bar\rho\|_2^2,
 \qquad
 \|P-\bar P\|_2^2+\|P-\bar P\|_4^4
                   \le C\|\rho-\bar\rho\|_2^2.
\]
Together with \eqref{uni:comparison}, this yields
\begin{equation}\label{lt:energy-tail}
 \mathcal H(t)+\|\rho(t)-\bar\rho\|_2^2
       +\int_t^\infty\bigl(D(s)+A_1^2(s)\bigr)\,\mathrm ds
       \le Ce^{-ct},\qquad t\ge0.
\end{equation}
Interpolation with the density upper bound gives
\begin{equation}\label{lt:density-decay}
 \|\rho(t)-\bar\rho\|_p
       \le C_pe^{-c\min\{1/2,1/p\}t},\qquad 1\le p<\infty.
\end{equation}

Since $E_*$ is uniformly bounded, \eqref{uni:logenergy-differential}
reduces to
\[
 (A_1^2)'+cy\le CD(t).
\]
For $t\ge1$, choose a time in $[t-1,t]$ at which $A_1^2$ does not
exceed its integral over that interval. Integrate from this time to
$t$, and then from $t$ to infinity. Equations
\eqref{lt:energy-tail} and \eqref{uni:comparison} give
\begin{equation}\label{lt:viscous-decay}
 A_1^2(t)+A_3^2(t)+\int_t^\infty y(s)\,\mathrm ds
       \le Ce^{-ct},\qquad t\ge1.
\end{equation}

For the corrected acceleration from \eqref{uni:qbounds},
$z\le2y+CA_3^4$, so \eqref{lt:viscous-decay} implies
$\int_t^\infty z(s)\,\mathrm ds\le Ce^{-ct}$.
The uniform acceleration bound gives $H\le C$ for $t\ge1$;
thus \eqref{uni:zdifferential} becomes
\[
 z'+cD_v\le CD(t),\qquad t\ge1.
\]
The same interval-averaging argument gives
\[
 z(t)+\int_t^\infty D_v(s)\,\mathrm ds\le Ce^{-ct},\qquad t\ge2.
\]
Now $y\le2z+CA_3^4$, and \eqref{uni:qbounds} gives
$\|\nabla q\|_2^2\le CD(t)$ for $t\ge1$.
Using $\dot u=v-q$, \eqref{uni:Hcomparison}, and the mass-weighted
Poincare inequality, we obtain
\begin{equation}\label{lt:force-decay}
 y(t)+\int_t^\infty\|\dot u(s)\|_{H^1}^2\,\mathrm ds
       \le Ce^{-ct},\qquad\hbox{for a.e. }t\ge1.
\end{equation}
Here the constant is enlarged on $[1,2]$ using \eqref{str:uniform}.

The low-exponent elliptic estimates imply
\[
 \|G\|_{H^1}+\|\omega\|_{H^1}
       \le C(A_1+A_2+A_3).
\]
For $p\ge2$, \eqref{mom:Hodge} and
$\diver u=(G+P-\bar P)/h$ therefore give
\[
 \|u\|_{W^{1,p}}
       \le C_p\bigl(A_1+A_2+A_3+\|\rho-\bar\rho\|_p\bigr).
\]
This proves \eqref{lt:relaxation}; for $p<2$ use the finite measure
of $\Om$. The positive-time $H^1$ continuity of $u$ and weak lower
semicontinuity in $W^{1,p}$ extend the estimate to every $t\ge1$.
Also $u_t=\dot u-u\cdot\nabla u$ and
\[
 \|\sqrt\rho\,(u\cdot\nabla u)\|_2
       \le C\|u\|_\infty\|\nabla u\|_2
\]
give \eqref{lt:acceleration}. In particular $u(t)$ tends to zero
uniformly. Since $|\bar P-\bar\rho^\gamma|\le C\|\rho-\bar\rho\|_2$,
the uncentered flux satisfies
$\|F(t)+\bar\rho^\gamma\|_{H^1}\le Ce^{-ct}$ for a.e. $t\ge1$.

We also need the time-integrated flux bounds. Sobolev embedding,
\eqref{mom:fluxelliptic}, and \eqref{lt:force-decay} imply, for each
fixed $2\le p<\infty$,
\begin{equation}\label{lt:flux-tail}
 \int_t^\infty\left(\|\rho\dot u\|_p^2
              +\|G\|_{W^{1,p}}^2+\|\omega\|_{W^{1,p}}^2\right)\,\mathrm ds
       \le C_pe^{-c_pt},\qquad t\ge1.
\end{equation}
Taking $p=4$ and summing over unit time intervals gives
\begin{equation}\label{lt:flux-integral}
\begin{aligned}
 \int_t^\infty(\|G\|_\infty+\|\omega\|_\infty)\,\mathrm ds
 &\le\sum_{j=0}^\infty
   \left[\int_{t+j}^{t+j+1}
       (\|G\|_\infty+\|\omega\|_\infty)^2\,\mathrm ds\right]^{1/2}\\
 &\le Ce^{-ct},\qquad t\ge1.
\end{aligned}
\end{equation}
Near the initial time, \eqref{eng:highgradient} and interpolation give
\[
 \|G\|_\infty
       \le C\bigl(1+y^{1/3}+\|\nabla\dot u\|_2^{2/3}\bigr).
\]
Weighted H\"older and \eqref{str:uniform}, as in Lemma~\ref{cl:initial},
show that $\int_0^1\|G\|_\infty\,\mathrm dt\le C$.

\subsection{Strictly positive initial density}

Assume $\min\rho_0>0$. On a trajectory where $\rho\le\bar\rho/2$,
Jensen's inequality and the continuity equation give
\[
 \Dt\theta(\rho)=\bar P-P-G\ge-\|G\|_\infty.
\]
Integrating from the initial time or the last crossing of
$\bar\rho/2$ yields
\[
 \theta(\rho(t,x))\ge
 \min\{\theta(\min\rho_0),\theta(\bar\rho/2)\}
                  -\int_0^\infty\|G\|_\infty\,\mathrm ds.
\]
Since $\theta(r)\to-\infty$ as $r\downarrow0$, the density has a
time-uniform positive lower bound.
By \eqref{eq:coefficient} and \eqref{eq:coeffconsequence},
\[
 |\lambda'(\rho)|\le C_\lambda\frac{1+\lambda(\rho)}\rho
       \le\frac{C}{\inf_{t,x}\rho(t,x)}.
\]
Thus all local coefficient bounds below are controlled by the stated
data. The finite-time bound for $\Phi(1)$ can also be restricted to
this positive density interval.
On its bounded positive range, monotonicity of $r^\gamma$ and
$\theta$ gives
\[
 \Dt|\theta(\rho)-\theta(\bar\rho)|
       +c|\theta(\rho)-\theta(\bar\rho)|
       \le|\bar P-\bar\rho^\gamma|+\|G\|_\infty.
\]
The mean pressure error decays exponentially. Solving this scalar
inequality and splitting the forcing integral at $t/2$, using
\eqref{lt:flux-integral}, proves
\begin{equation}\label{lt:density-sup}
 \|\rho(t)-\bar\rho\|_\infty\le Ce^{-ct}.
\end{equation}

To control density gradients, retain the pressure term in
\eqref{str:Phi}. Since $\gamma\rho^\gamma/h$ is bounded below,
\begin{equation}\label{lt:gradient-damping}
 \frac{\mathrm d}{\mathrm dt}\|\Phi\|_q+c\|\Phi\|_q
       \le C\|\nabla u\|_\infty\|\Phi\|_q
                      +C\|\rho\dot u\|_q+Ce^{-ct}.
\end{equation}
The flux estimates and \eqref{lt:relaxation} give
\[
 \|\nabla G\|_q+\|\nabla\omega\|_q+\|G\|_\infty+\|\omega\|_\infty
       \le C_q(\|\rho\dot u\|_q+e^{-ct}).
\]
Differentiating $\diver u=(G+P-\bar P)/h$ once therefore yields
\[
 \|u\|_{W^{2,q}}
       \le C_q\bigl[(1+\|\rho\dot u\|_q)\|\Phi\|_q
                         +\|\rho\dot u\|_q+e^{-ct}\bigr].
\]
Apply \eqref{str:BKM} to $e^tu(t)$. Together with
\eqref{lt:density-sup} and the preceding estimate, this gives
\begin{equation}\label{lt:integrable-lip}
\begin{aligned}
 \|\nabla u\|_\infty
 &\le b(t)\bigl[1+\log(e+\|\Phi\|_q)\bigr],\\
 b(t)&=C_q(1+t)\bigl(e^{-ct}+\|\rho\dot u\|_q
                                      +\|\rho\dot u\|_q^2\bigr),\\
 \int_t^\infty b(s)\,\mathrm ds&\le C_qe^{-ct},\qquad t\ge1.
\end{aligned}
\end{equation}
The last bound follows from \eqref{lt:flux-tail} by the same
unit-interval summation. All coefficients here are controlled
independently of $\|\Phi\|_q$.
Divide \eqref{lt:gradient-damping} by $e+\|\Phi\|_q$ and first
discard its damping term. Logarithmic Gronwall gives
$\sup_{t\ge1}\|\Phi(t)\|_q<\infty$.
Hence \eqref{lt:integrable-lip} gives an exponential $L^1$ tail for
$\|\nabla u\|_\infty$. Returning to
\eqref{lt:gradient-damping} and using \eqref{lt:flux-tail} now gives
\[
 \|\nabla\rho(t)\|_q\le Ce^{-ct}.
\]
This proves the density assertion in \eqref{lt:positive}.

The $H^2$ div--curl estimate gives
\[
\begin{aligned}
 \|u\|_{H^2}
 &\le C\bigl(A_2+\|u\|_{H^1}+\|\nabla\rho\|_2\bigr)
       +C\|\nabla\rho\|_q\|\diver u\|_{2q/(q-2)},\\
 \|\diver u\|_{2q/(q-2)}
 &\le C\|u\|_{H^1}^{1-2/q}\|u\|_{H^2}^{2/q}+C\|u\|_{H^1}.
\end{aligned}
\]
Young's inequality absorbs the $H^2$ term, proving
$\|u(t)\|_{H^2}\le Ce^{-ct}$.
The equations $u_t=\dot u-u\cdot\nabla u$ and
$\rho_t=-u\cdot\nabla\rho-\rho\diver u$ give the remaining pointwise
bounds in \eqref{lt:positive}.
The preceding $W^{2,q}$ estimate, \eqref{lt:flux-tail}, and
\[
 \|u\cdot\nabla u\|_{H^1}\le C\|u\|_{H^2}^2
\]
give its time-integrated bounds. This proves \eqref{lt:positive}.

\subsection{Persistence of vacuum and density gradients}

For each finite $T$, \eqref{str:densitygradient} gives a Lipschitz
flow preserving $\overline\Om$. If $X(t,a)$ is the trajectory from
$a\in\overline\Om$, then
\[
 \rho(t,X(t,a))=\rho_0(a)
       \exp\left(-\int_0^t\diver u(s,X(s,a))\,\mathrm ds\right).
\]
Thus an initial vacuum point persists, proving \eqref{lt:vacuum}.
For $2<r\le q$, the mean-zero Gagliardo--Nirenberg inequality gives
\[
 \bar\rho\le\|\rho-\bar\rho\|_\infty
       \le C_r\|\nabla\rho\|_r^{r/[2(r-1)]}
                   \|\rho-\bar\rho\|_2^{(r-2)/[2(r-1)]}.
\]
Together with \eqref{lt:energy-tail}, this proves
\eqref{lt:gradient-growth}. Finally,
\[
 |\{\rho(t)\le\varepsilon\bar\rho\}|
       \le\frac{\|\rho(t)-\bar\rho\|_2^2}
                    {(1-\varepsilon)^2\bar\rho^2}
       \le C_\varepsilon e^{-ct},\qquad 0\le\varepsilon<1,
\]
which proves \eqref{lt:vacuum-set} and completes the proof of
Theorem~\ref{thm:longtime}.

\section*{Acknowledgments}
The authors acknowledge the use of generative AI tools in developing and
checking mathematical arguments and in drafting and revising the manuscript.

\end{document}